\documentclass[12pt]{amsart}
\usepackage{amsmath}
\usepackage{amscd}
\usepackage{amssymb}
\usepackage{pdfsync}

\newtheorem{theorem}{Theorem}[section]

\newtheorem{proposition}[theorem]{Proposition}

\newtheorem*{lem}{Lemma}

\theoremstyle{definition}

\newcommand{\cO}{{\mathcal O}}

\newcommand{\cN}{{\mathcal N}}

\newcommand{\tr}{{\text{tr}}}

\newcommand{\Hom}{{\text{Hom}}}

\newcommand{\End}{{\text{End}}}
\newcommand{\Ad}{{\text{Ad}}}

\newcommand{\Ima}{{\text{Im}}}
\newcommand{\spn}{{\text{span}}}

\newcommand{\Lb}{{\mathfrak{b}}}

\newcommand{\Lg}{{\mathfrak g}}
\newcommand{\Lo}{{\mathfrak o}}

\newcommand{\LS}{{\mathfrak{S}}}

\newcommand{\bV}{{\mathbf{V}}}

\newcommand{\bW}{{\mathbf{W}}}

\newcommand{\tF}{{\textbf{F}}}
\newcommand{\tk}{{\textbf{k}}}

\newcommand{\p}{\perp}

\newcommand{\bbeta}{\bar{\beta}}
\newcommand{\bbW}{\bar{\bW}}
\newcommand{\tbV}{\tilde{\bV}}
\newcommand{\brA}{\bar{A}}
\newcommand{\brm}{\bar{m}}
\newcommand{\bv}{\bar{v}}
\newcommand{\bu}{{\bar{u}}}
\newcommand{\beq}{\begin{equation*}}
\newcommand{\eeq}{\end{equation*}}

\begin{document}
\title[Nilpotent pieces in the dual of odd orthogonal Lie algebras]{Nilpotent pieces in the dual of odd orthogonal Lie algebras}
        \author{Ting Xue}
        \address{Department of Mathematics, Northwestern University,
Evanston, IL 60208, USA}
        \email{txue@math.northwestern.edu}

\maketitle

\begin{centerline}{\em To George Lusztig on the occasion of his 65th birthday}
\end{centerline}

\begin{abstract}
Let $\cN_{\Lg^*}$ be the variety of nilpotent elements in the dual of the Lie algebra of a reductive algebraic group over an algebraically closed field. In \cite{Lu2} Lusztig proposes a definition of a  partition of $\cN_{\Lg^*}$ into smooth locally closed  subvarieties (which are indexed by the unipotent classes in the corresponding group over complex numbers) and gives explicit results in types $A$, $C$ and $D$. We discuss type $B$ in this note.
\end{abstract}

\section{Introduction}
Let $\tk$ be an algebraically closed field of characteristic $p\geq 0$. Let $G$ be a connected reductive algebraic group over $\tk$ and $\Lg$ the Lie algebra of $G$. Let $\Lg^*$ be the dual vector space of $\Lg$. Note that $G$ acts on ${\Lg^*}$ by coadjoint action. Let $\cN_{\Lg^*}$ be the variety of nilpotent elements in $\Lg^*$ (where an element $\xi:\Lg\rightarrow\tk$ is called nilpotent if it annihilates some Borel subalgebra of $\Lg$, see \cite{KW}).   Let $G_{\mathbb{C}}$ be the reductive group over $\mathbb{C}$ of the same type as $G$. In \cite{Lu2} Lusztig proposes a definition of a  partition of $\cN_{\Lg^*}$ into smooth locally closed $G$-stable pieces, called nilpotent pieces which are indexed by the unipotent classes in $G_{\mathbb{C}}$. The case where $G$ is of type $A$, $C$, or $D$ has been illustrated in \cite{Lu2}. We treat in this note the case where $G$ is of type $B$.

In section 2 we recall the definition of nilpotent pieces in $\cN_{\Lg^*}$ that Lusztig proposes and state the main theorem. We also include Lusztig's result (and proof) on the number of rational points in $\cN_{\Lg^*}$ over finite fields of characteristic 2 (where $G$ is of type $B$, see \ref{sec-number}). In section 3 we give an explicit description of the set $\Lg_2^{*\delta!}$ (see 2.1) involved in the definition of nilpotent pieces, where a key definition is suggested by Lusztig (see \ref{sec-def}). The main theorem is proved in section 4. We give some examples in section 5.

After this paper has been submitted, a preprint \cite{CP} of Clarke and Premet appears which shows in a uniform way  that the definition of unipotent pieces  in $G$ proposed by Lusztig in \cite{Lu4} works for all types and which proves analogous results for nilpotent pieces in $\Lg$ and $\Lg^*$ (for all types). Our approach is based on another definition  proposed by Lusztig in \cite{Lu2} and gives explicit description of the nilpotent pieces in $\mathfrak{so}({2n+1})^*$ which implies the smoothness of nilpotent pieces in this case. One can also describe explicitly
which nilpotent orbits in $\Lg^*$ (for type $B,C$) lie in the same piece using the construction of \cite{Lu2} and this paper (this will be done elsewhere).   It has been  pointed out in \cite{CP} that the definitions proposed by Lusztig in \cite{Lu4,Lu2} give rise to the same partition of $\cN_{\Lg^*}$ into pieces for classical groups in view of \cite{CP,Lu2} and our result.

\vskip 10pt {\noindent\bf\large Acknowledgement}  \ I wish to thank George Lusztig for suggesting a key definition, for allowing me to include  his proof on the number of nilpotent elements in $\Lg^*$, and for many helpful discussions. I am also grateful to the referees for suggestions and comments.
\section{Preliminaries and statement of the main theorem}
\subsection{}\label{ssec-d1} In this subsection we recall the definition of nilpotent pieces in $\cN_{\Lg^*}$ that Lusztig proposes (see \cite{Lu2} for more details). Let  $\mathfrak{D}_{G_{\mathbb{C}}}$ be the set of all $f\in\text{Hom}(\mathbb{C}^*,G_{\mathbb{C}})$ such that there exists a homomorphism of algebraic groups $\tilde{f}:SL_2(\mathbb{C})\to G_{\mathbb{C}}$ with $\tilde{f}\left(\begin{array}{cc}r&0\\0&r^{-1}\end{array}\right)=f(r)$ for all $r\in\mathbb{C}^*$. Let $\mathfrak{D}_G$ be the set of all $\delta\in\Hom(\tk^*,G)$ such that the image of $\delta$ in $G\backslash\Hom(\tk^*,G)=G_{\mathbb{C}}\backslash\Hom(\mathbb{C}^*,G_{\mathbb{C}})$ can be represented by an element in $\mathfrak{D}_{G_{\mathbb{C}}}$.

Let $\delta\in\mathfrak{D}_G$. We have $\Lg=\oplus_{i\in\mathbb{Z}}\Lg_i^\delta$, where $\Lg_i^\delta=\{x\in\Lg|\Ad(\delta(r))x=r^ix,\ \forall\ r\in\tk^*\}$. We set $\Lg_{\geq i}^\delta=\oplus_{i'\geq i}\Lg_{i'}^\delta$.  For $j\in\mathbb{Z}$, let $\Lg_j^{*\delta}=\text{Ann}(\oplus_{i\neq -j}\Lg^\delta_i)\subset\Lg^*$. We have $\Lg^*=\oplus_{j\in\mathbb{Z}}\Lg_j^{*\delta}$.  For any $\xi\in\Lg^*$, denote $Z_G(\xi)$ the centralizer of $\xi$ in $G$ under the coadjoint action. Let
\begin{equation*}\label{eqn-g2d}\Lg_2^{*\delta!}=\{\xi\in\Lg_2^{*\delta}|Z_G(\xi)\subset G_{\geq 0}^\delta\},
\end{equation*}
where $G_{\geq 0}^\delta$ is the (well-defined) closed connected subgroup of $G$ such that its Lie algebra is $\Lg_{\geq 0}^\delta$.

Let $D_G$ be the set of equivalence classes in $\mathfrak{D}_G$, where $\delta$ and $\delta' $ in $\mathfrak{D}_G$ are said to be equivalent if for any $i\in\mathbb{N}=\{0,1,2,\ldots\}$,  $\Lg_{\geq i}^\delta=\Lg_{\geq i}^{\delta'}$. Let $\Delta\in D_G$. We write $G_{\geq 0}^\Delta$, $\Lg_{\geq i}^\Delta$ instead of $G_{\geq 0}^\delta$, $\Lg_{\geq i}^\delta$ for $\delta\in\Delta$. For $j\in\mathbb{N}$, let $\Lg_{\geq j}^{*\Delta}=\text{Ann}(\Lg^\Delta_{\geq -j+1})$. Then $\Lg_{\geq j}^{*\Delta}=\oplus_{j'\geq j}\Lg_{j'}^{*\delta}$ ($\delta\in\Delta$). For any $\delta\in\Delta$, let $\Sigma^{*\delta}\subset\Lg^{*\Delta}_{\geq 2}/\Lg^{*\Delta}_{\geq 3}$ be the image of $\Lg_2^{*\delta!}\subset \Lg_2^{*\delta}$ under the obvious isomorphism $\Lg_2^{*\delta}\xrightarrow{\sim}\Lg^{*\Delta}_{\geq 2}/\Lg^{*\Delta}_{\geq 3}$. Then $\Sigma^{*\delta}$ is independent of the choice of $\delta$ in $\Delta$; we denote it  by $\Sigma^{*\Delta}$. Let $\sigma^{*\Delta}\subset\Lg_{\geq 2}^{*\Delta}$ be the inverse image of $\Sigma^{*\Delta}$ under the obvious map $\Lg_{\geq 2}^{*\Delta}\to\Lg^{*\Delta}_{\geq 2}/\Lg^{*\Delta}_{\geq 3}$. Then $\sigma^{*\Delta}$ is stable under the coadjoint action of $G_{\geq 0}^\Delta$ on $\Lg^{*\Delta}_{\geq 2}$. We have a map
$$\Psi_{\Lg^*}:\sqcup_{\Delta\in D_G}\sigma^{*\Delta}\rightarrow\cN_{\Lg^*},\ \xi\mapsto\xi.$$

\begin{theorem}\label{thm}If $G$ is of type $B$, then the map
$\Psi_{\Lg^*}$
is a bijection.
\end{theorem}

In \cite{Lu2} Lusztig conjectures that $\Psi_{\Lg^*}$ is a bijection for any $G$ and proves this in the case where $G$ is of type $A$, $C$ or $D$.
Theorem \ref{thm} will be proved in section \ref{sec-thm}. We can assume that $G$ is an odd special orthogonal group.

Let $\mathfrak{U}_G$ be the set of $G$-orbits on $D_G$. Then $\mathfrak{U}_G$ is a finite set that depends only on the type of $G$, not on $\tk$ (see \cite{Lu1}). For any $\cO\in\mathfrak{U}_G$, we set
$$\cN_{\Lg^*}^\cO=\Psi_{\Lg^*}(\sqcup_{\Delta\in\cO}\sigma^{*\Delta}).$$
The subsets $\cN_{\Lg^*}^\cO$ are called {\em pieces} of $\cN_{\Lg^*}$. They form a partition of $\cN_{\Lg^*}$ into smooth  locally closed subvarieties (which are unions of $G$-orbits) indexed by $\mathfrak{U}_G=\mathfrak{U}_{G_\mathbb{C}}$.

\subsection{}\label{ssec-v} Let $\bV$ be a
vector space   over $\tk$ equipped with a  nondegenerate quadratic form $Q:\bV\rightarrow \tk$. Let $\beta:\bV\times \bV\to\tk $ be the bilinear
form associated to $Q$, namely, $\beta(v,v') =Q(v+v')-Q(v)-Q(v')$ for
$v,v'\in \bV$.  Then $Q:\text{Rad}(Q)=\{v\in \bV|\beta( v,\bV)  =0\}\rightarrow\tk$ is
injective and $\text{Rad}(Q)=0$ unless $p=2$ and $\dim\bV$ is odd.
The special orthogonal group $SO(\bV)$ is the identity component of $O(\bV)=\{g\in\text{GL}(\bV)|\ Q(gv)=Q(v),\ \forall\ v\in \bV\}$ and its Lie algebra is $\mathfrak{o}(\bV)=\{x\in\End (\bV)|\ \beta(xv,v)=0,\
\forall\ v\in \bV\text{ and }x|_{\text{Rad}(Q)}=0\}$.

A $\mathbb{Z}$-grading $\bV=\oplus_{a\in\mathbb{Z}}\bV^a$ of $\bV$ is called an {\em o-good} grading if $\dim \bV^a=\dim \bV^{-a}\geq\dim \bV^{-a-2}$ for all $a\geq 0$, $\dim \bV^a$ is even for all odd $a$, $\beta( \bV^a,\bV^b)=0$
whenever $a+b\neq 0$, and $Q|_{\bV^a}=0$ for all $a\neq 0$.

A filtration  $\bV_*=(\bV^{\geq a})_{a\in\mathbb{Z}}$ of $\bV$ (where  $\bV^{\geq a+1}\subset\bV^{\geq a}$, $\bV^{\geq a}=\{0\}$ for some $a$ and $\bV^{\geq a}=\bV$ for some $a$) is called a $Q$-filtration
 if  $Q|_{\bV^{\geq a}}=0$ and
$\bV^{\geq 1-a}=(\bV^{\geq a})^\p$ for any $a\geq 1$.

\subsection{}\label{sec-ogood}

From now on we assume that $\bV$ is of dimension $2N+1$, $Q$ is a fixed nondegenerate quadratic form on $\bV$ with associated bilinear form $\beta$ and  that $G=SO(\bV)$ (with respect to $Q$). Let $R$ denote  the radical $\text{Rad}(Q)$ of $Q$. For any subspace $\bW\subset\bV$, let $\bW^\p$ denote the set $\{v\in\bV|\beta(v,\bW)=0\}$.

To give an element $\delta \in\mathfrak{D}_G$ is the same as to give an $o$-good grading $\bV=\oplus_{a\in\mathbb{Z}}\bV^a$ of $\bV$ ($\delta$ is given by $\delta(r)|_{\bV^a}=r^a$ for all $r\in\tk^*$ and all $a\in\mathbb{Z}$, see \cite[1.5]{Lu1}). Let $\mathfrak{F}_o(\bV)$ be the set of all $Q$-filtrations $\bV_*=(\bV^{\geq a})_{a\in\mathbb{Z}}$ such that there exists an $o$-good grading $\oplus_{a\in\mathbb{Z}}\bV^a$ of $\bV$ with $\bV^{\geq a}=\oplus_{a'\geq a}\bV^{a'}$.
The set $\mathfrak{F}_o(\bV)$ and the set $D_G$ (see \ref{ssec-d1}) are identified as follows  (see \cite[2.7]{Lu1}). Let $[\delta]\in D_G$ be the equivalence class containing $\delta\in\mathfrak{D}_G$ and let $\bV=\oplus_{a\in\mathbb{Z}}\bV^a$ be the $o$-good grading corresponding to $\delta$. Then $\bV_*=(\bV^{\geq a})$ with $\bV^{\geq a}=\oplus_{a'\geq a}\bV^{a'}$ is the element of $\mathfrak{F}_o(\bV)$ corresponding to $[\delta]$.

\subsection{}\label{ssec-map}
Let $\LS(\bV)$ denote the set of all symplectic bilinear forms on $\bV$. In \cite[3.1]{X2}, we have defined a map (assume $p=2$)
$$\Lg^*\to\LS(\bV),\ \xi\mapsto\beta_\xi,$$
where $\beta_\xi(v,v')=\beta(Xv,v')-\beta(v,Xv')$, $\forall\ v,v'\in\bV$ and $X$ is such that $\xi(x)=\text{tr}(Xx)$, $\forall\ x\in\Lg$. In fact for arbitrary $p$ this map makes sense and is a vector space isomorphism. In the following  we denote by $\beta_\xi$ the symplectic bilinear form that corresponds to $\xi\in\Lg^*$ under this map.

Let $\LS(\bV)_{nil}$ denote the set of all symplectic bilinear forms $\beta_\xi$ with $\xi\in\cN_{\Lg^*}$.

Assume $p=2$ and $\beta_\xi\in \LS(\bV)_{nil}$. There are  a unique $m\in\mathbb{N}$ and a unique set of vectors  $\{v_i,i\in[0,m]\}$  such that (see \cite[Lemma 3.5]{X2})
 \begin{eqnarray*}
 &&\mathrm{(a)}\quad\beta(v_m,v)=0,\ \beta_\xi(v_i,v)=\beta(v_{i-1},v), i\in[1,m],\ \beta_\xi{(v_0,v)}=0, \forall\ v\in \bV;\\
 &&\qquad Q(v_i)=0, i\in[0,m-1],\  Q(v_m)=1.
 \end{eqnarray*}

\subsection{}\label{sec-number}The proofs in this subsection are due to G. Lusztig.

A basis $(e_i)_{i\in[-N,N]}$ of $\bV$ is said to be {\em good} if $\beta(e_i,e_j)=\delta_{i+j,0}+\delta_{i,0}\delta_{j,0}$ for all $i,j\in[-N,N]$, and $Q(e_i)=\delta_{i,0}$ for all $i\in[-N,N]$.

We show that the following two conditions  are equivalent:\\[5pt]
\indent\quad(a1) $\beta_\xi\in\LS(\bV)_{nil}$;\\[5pt]
\indent\quad(a2) there exists a {\em good} basis $(e_i)_{i\in[-N,N]}$ of $\bV$ such that $i+j\geq 0$ implies $\beta_\xi(e_i,e_j)=0$.\\[5pt]
Assume $\beta_\xi\in\LS(\bV)_{nil}$. We can find a Borel subalgebra $\Lb$ of $\Lg$ such that $\xi(\Lb)=0$. There exists a {\em good} basis $(e_i)_{i\in[-N,N]}$ of $\bV$ such that for all $x\in\Lb$, $xe_i=\sum_{j\geq i}x_{ij}e_j$ for all $i\in[-N,N]$, where $x_{ij}+x_{-j,-i}=0$ for all $i,j\in[-N,N]-\{0\}$, $x_{i,-i}=0$ for all $i\in[-N,N]-\{0\}$,  $x_{0,i}+2x_{-i,0}=0$ for all $i\in[-N,N]-\{0\}$ and $2x_{0,0}=0$. We take $X\in\End(\bV)$ such that $\xi(x)=\tr(Xx)$ for all $x\in\Lg$. Assume $Xe_i=\sum_jX_{ij}e_j$. For $x\in\Lb$, we have $\tr(Xx)=\sum_{-N\leq j\leq i\leq -1}(X_{ij}-X_{-j,-i})x_{ji}+\sum_{-N\leq j<i\leq -1}(X_{-i,j}-X_{-j,i})x_{j,-i}+\sum_{j<0}(X_{0,j}-2X_{-j,0})x_{j,0}$.  It follows from $\xi(\Lb)=0$ that $X_{ij}-X_{-j,-i}=0$ for  all $i\geq j$, and $X_{0j}-2X_{-j,0}=0$ for all $j\in[-N,-1]$. Now we have $\beta_\xi(e_i,e_j)=X_{i,-j}-X_{j,-i}=0$ for $i,j\in[-N,N]-\{0\}$, $i\geq -j$; $\beta_\xi(e_0,e_0)=0$ and $\beta_\xi(e_0,e_j) =X_{0,-j}-2X_{j0}=0$ for $j\in[1,N]$.
Thus (a2) holds. Conversely assume (a2) holds. Let $(e_i)_{i\in[-N,N]}$ be a basis of $\bV$ as in (a2). Let $X$ and $X_{ij}$ be as in the first part of the proof. We have for $i,j\in[-N,N]-\{0\}$, $i\geq -j$, $X_{i,-j}-X_{j,-i}=0$; and for $j\in[1,N]$,  $X_{0,-j}-2X_{j0}=0$. Let $\Lb$ be a Borel subalgebra related to $(e_i)$ as in the first part of the proof. Then $\xi(\Lb)=0$. Hence (a1) holds.

We assume now that $\tk$ is an algebraic closure of the finite prime field $\tF_2$ and that $\dim\bV\geq 3$. We choose an $\tF_2$ rational structure on $\bV$ such that $Q$ is defined over $\tF_2$. Then the Frobenius map $F$ relative to this $\tF_2$ structure acts naturally and compatibly on $\LS(\bV)_{nil}$. We show in the remainder of this subsection that

\quad (b) {\em $|\LS(\bV)^{F^n}_{nil}|=q^{2N^2}, \text{ where } q=2^n.$}

For any $m\in[0,N]$, let $S_m$ be the set of all sequences $v_*=(v_0,\ldots,v_{m-1})$ of linearly independent vectors in $\bV$ such that $Q|_{\spn\{v_i,i\in[0,m-1]\}}=0$.  We have

\quad(c) {\em $|S_m^{F^n}|=(q^{2N}-1)(q^{2N-2}-1)\cdots(q^{2N-2m+2}-1)q^{m(m-1)/2}.$}\\[10pt]
For each $v_*\in S_m$, let $\cN_{v_*}$ be the set of nilpotent elements in $\Lo(\bV')=\{T\in\End(\bV')|\beta'(Tv',v')=0,\ \forall \ v'\in\bV'\}$, where $\bV'=L^\p/L$ with $L=\spn\{v_i,i\in[0,m-1]\}\oplus R$ and $\beta'$ is the bilinear form on $\bV'$ induced by $\beta$.
By a result of Springer (see \cite{Sp}), we have

\quad (d) {\em $
|\cN_{v_*}^{F^n}|=q^{2(N-m)(N-m-1)}.$}

Fix a pair $(v_*,T)$ where $v_*\in S_m$ and $T\in\cN_{v_*}$.   Let $\mathcal{E}_{(v_*,T)}$ be the set of all  $(,)\in\LS(\bV)$ such that $(v_0,v)=0$, $(v_i,v)=\beta(v_{i-1},v)$, $i\in[1,m-1]$, for all $v\in\bV$;
 $(r,v)=\beta(v_{m-1},v)$ for all $v\in\bV$, where $r\in R$ is such that $Q(r)=1$;
$\beta'(T v',v'')=(v',v'')'$,  $\forall\ v',v''\in\bV'$, where $(,)'$ is the bilinear form on $\bV'$ induced by  $(,)$ (note that $(L,L^\p)=0$).We show that\\[5pt]
\indent\quad (e) {\em $\mathcal{E}_{(v_*,T)}\subset\LS(\bV)_{nil}$;}\\[5pt]
\indent\quad (f) {\em $\mathcal{E}_{(v_*,T)}\text{ is an affine space of dimension }2m(N-m)+m(m-1)/2.$}\\[5pt]
Let $(,)\in\mathcal{E}_{(v_*,T)}$. Let $\tilde{\bV}$ be a complement of $L$ in $L^\p$ and let $L'$ be a complement to $L^\p$ in $\bV$ such that $Q|_{L'}=0$ and $\beta(L',\tilde{\bV})=0$. Note that $\beta|_{\tilde{\bV}}$ is nondegenerate and can be identified with $\beta'$ on $\bV'$. We can regard $T\in\cN_{v_*}$ as an element $\tilde{T}\in\cN_{\Lo(\tilde{\bV})}$ and identify $(,)|_{\tilde{\bV}}$ with $(,)'$ on $\bV'$. Let $N'=N-m$. We can find a basis $(e_i)_{i\in[-N',N']-\{0\}}$ of $\tilde{\bV}$ such that $\beta(e_i,e_j)=\delta_{i+j,0}$ for all $i,j\in[-N',N']-\{0\}$, $Q(e_i)=0$ for all $i\in[-N',N']$, and $(e_i,e_j)=0$ for all $i+j\geq 0$ (since $(v,v')=\beta(\tilde{T}v,v')$ on $\tilde{\bV}$ and $\tilde{T}$ is nilpotent). We extend it to a basis of $\bV$ by setting $e_{N'+j}=v_{m-j}$, $j\in[1,m]$, $e_0=r$ and $e_{-N'-j}=u_{m-j}$, $j\in[1,m]$, where $u_{m-j}$ are elements in $L'$ such that $\beta(v_{m-j},u_{m-j'})=\delta_{j,j'}$. It follows that $(e_i)_{i\in[-N,N]}$ is a good basis of $\bV$. We show that $(e_i,e_j)=0$ if $i,j\in[-N,N]$ and $i+j\geq 0$. For $j,j'\in[1,m-1] $ and $j'\geq j$,
 $(e_{-N'-j},e_{N'+j'})=(u_{m-j},v_{m-j'})=\beta(u_{m-j},v_{m-j'-1})=0$; for $i\in[-N',N']-\{0\}$, $j\in[1,m-1]$, $(e_i,e_{N'+j})=(e_i,v_{m-j})=\beta(e_i,v_{m-j-1})=0$; for $j,j'\in[1,m-1]$, $(e_{N'+j},e_{N'+j'})=\beta(v_{m-j},v_{m-j'-1})=0$; for $j\in[0,N]$, $(e_0,e_j)=\beta(v_{m-1},e_j)=0$ and for $i\in[-N,N]$, $(e_i,e_{N'+m})=(e_i,v_{0})=0$. This completes the proof of (e) (we use (a1)$\Leftrightarrow$(a2)). We prove (f). Let $\tilde{\bV},L'$ be as in the proof of (e). Let $u_0,\ldots,u_{m-1}$ be a basis of $L'$  and let $(e_i)_{i\in[-N',N']-\{0\}}$ be a basis of $\tilde{\bV}$. To specify an element $(,)$ of $\mathcal{E}_{(v_*,T)}$, we need to specify $(u_i,e_j)$, $i\in[0,m-1]$, $j\in[-N',N']-\{0\}$ and $(u_i,u_j)$ for $i<j$ in $[0,m-1]$. These provide coordinates in $\mathcal{E}_{(v_*,T)}$. Thus (f) follows.

Let $\beta_\xi\in\LS(\bV)_{nil}$. Let $\{v_i,i\in[0,m]\}$ be the set of vectors as in \ref{ssec-map} (a). Then $v_i$, $i\in[0,m]$ are linearly independent, $v_m\in R$ and $Q|_{\spn\{v_i,i\in[0,m-1]\}}=0$ (see \cite[Lemmas 3.5, 3.7]{X2}). We set  $L=\spn\{v_i,i\in[0,m]\}$ and $\bV'=L^\p/L$. The bilinear form $\beta$ induces a nondegenerate bilinear form $\beta'$ on $\bV'$. For any $x\in L$,  $\beta_\xi(x,L^\p)=0$ since $\beta_\xi(v_0,\bV)=0$ and $\beta_\xi(v_i,L^\p)=\beta(v_{i-1},L^\p)=0$ for all $i\in[1,m]$. Hence $\beta_\xi$ induces a symplectic bilinear form $\beta_\xi'$ on $\bV'$. We define $T_\xi':\bV'\to\bV'$ by $\beta'(T_\xi' v',v'')=\beta_\xi'(v',v'')$ for any $v',v''\in\bV'$. Note that $T_\xi'\in\Lo(\bV')$ and $T_\xi'$ is nilpotent (see \cite[Lemma 3.11]{X2}). Thus we have a natural map $\beta_\xi\mapsto(v_*=(v_0,\ldots,v_{m-1}),T_\xi')$ from $\LS(\bV)_{nil}$ to the set of all pairs $(v_*,T)$, where $v_*\in S_m$ for some $m\in[0,N]$ and $T\in\cN_{v_*}$. The fiber of this map at $(v_*,T)$ is the set $\mathcal{E}_{(v_*,T)}$ (see (e)). Now using (c), (d) and (f) we get
\begin{eqnarray*}
&&|\LS(\bV)^{F^n}_{nil}|=\sum_{m\in[0,N]}\sum_{v_*\in S_m^{F^n}}q^{2(N-m)(N-m-1)}q^{2m(N-m)+m(m-1)/2}\\&&=q^{2N^2}\sum_{m\in[0,N]}(1-q^{-2N})(1-q^{-2N+2})\cdots(1-q^{-2N+2m-2})q^{-2N+2m}.
\end{eqnarray*}
Let $X_N=\sum_{m=0}^N(1-q^{N})(1-q^{N-1})\cdots(1-q^{N-m+1})q^{N-m}$. We have $X_1=1$ and $X_{N+1}=q^{N+1}+(1-q^{N+1})X_N$. Hence $X_N=1$ by induction on $N$ and (b) follows.

\section{The set $\Lg_2^{*\delta!}$}
\noindent In this section we fix $\delta\in\mathfrak{D}_G$ and the corresponding $o$-good grading $\bV=\oplus_{a\in\mathbb{Z}}\bV^a$ of $\bV$ (see \ref{sec-ogood}) and describe the set $\Lg_2^{*\delta!}$ (see \ref{ssec-d1}) more explicitly (see Proposition \ref{prop-2}). 
\subsection{}\label{sec-def}
Let $\mathfrak{S}(\bV)_2$ be the set of all symplectic bilinear forms $\beta_\xi\in\LS(\bV)$ such that $\beta_\xi(\bV^a,\bV^b)=0$, whenever $a+b\neq -2$.  For $\beta_\xi\in\LS(\bV)_2$, define $A:\bV^a\rightarrow \bV^{a+2}$ for all $a\neq -2$ by
$$\beta(Ax^a,x^{-a-2})=\beta_\xi(x^a,x^{-a-2}), \forall\ x^a\in \bV^a,x^{-a-2}\in \bV^{-a-2}.$$
Let $\LS(\bV)_2^0$ be the set of all $\beta_\xi\in\LS(\bV)_2$ such that the following conditions (a) and (b) hold  (this definition is suggested by Lusztig):

(a) all maps $\bV^0\xrightarrow {A}\bV^2\xrightarrow{A}\cdots\xrightarrow{A}\bV^{2n}
\xrightarrow{A}\cdots$ are surjective, and $Q|_{\ker(A^n:\bV^0\rightarrow \bV^{2n})}$ is nondegenerate for all $n\geq 1$;\\[5pt]
\indent(b) all maps $\bV^{-1}\xrightarrow{A}\bV^1\xrightarrow{A}\cdots\xrightarrow{A}\bV^{2n-1}\xrightarrow{A}\cdots$ are surjective, and   the symplectic form $(-,-)_{-1}:=\beta(A-,-)$ on $\ker(A^n:\bV^{-1}\rightarrow \bV^{2n-1})$ is nondegenerate for all $n\geq 1$.

\begin{proposition}\label{prop-2}We have $\xi\in\Lg_2^{*\delta !}\text{ if and only if }\beta_\xi\in\mathfrak{S}(\bV)_2^0.$\end{proposition}

Note that $\xi\in\Lg_2^{*\delta}$ if and only if $\beta_\xi\in\mathfrak{S}(\bV)_2$ (since $\delta(r)|_{\bV^a}=r^a$, and $\xi\in\Lg_2^{*\delta}$ iff  $\beta_\xi(\delta(r)v,\delta(r)v')=r^{-2}\beta_\xi(v,v')$ for all $v,v'\in\bV$,  $r\in\tk^*$). The proof of the proposition in the case when $p\neq 2$ is  given in subsection \ref{sec-charnot2} and that in the case when $p=2$ is given in subsections \ref{sec-1}-\ref{sec-2}. In general, $x^k$ denotes an element in $\bV^k$.

We first show that the condition (b) can be reformulated as follows

(b$'$)  $A^{2n-1}: \bV^{-2n+1}\rightarrow \bV^{2n-1}$ is an isomorphism for all $n\geq 1$.\\[10pt]
Let
$K_{2n-1,\bV}=\ker(A^n:\bV^{-1}\rightarrow \bV^{2n-1})$ and $I_{2n-1,\bV}=\Ima(A^{n-1}:\bV^{-2n+1}\to\bV^{-1})$.
Assume (b) holds. We show that $A^{2n-1}:\bV^{-2n+1}\rightarrow \bV^{2n-1}$ is injective. Assume $A^{2n-1}x^{-2n+1}=0$. Then $A^{n-1}x^{-2n+1}\in\text{Rad}((,)_{-1}|_{K_{2n-1,\bV}})=\{0\}$ and  then $\beta(x^{-2n+1},\bV^{2n-1})=\beta(x^{-2n+1},A^n\bV^{-1})=\beta(A^nx^{-2n+1},\bV^{-1})=0$ (note that $A^n:\bV^{-1}\to\bV^{2n-1}$ is surjective). Hence $x^{-2n+1}=0$ (since $\beta:\bV^a\times\bV^{-a}\to\tk$ is nondegenerate for all $a\neq 0$). Now (b$'$) follows since $\dim\bV^{-2n+1}=\dim \bV^{2n-1}$.

Conversely assume (b$'$) holds. It is clear that all maps $A$ in (b) are surjective.  We have $\bV^{-1}=I_{2n-1,\bV}\oplus K_{2n-1,\bV}$ (since $I_{2n-1,\bV}\cap K_{2n-1,\bV}=\{0\}$,  $\dim I_{2n-1,\bV}=\dim\bV^{-2n+1}$ and $\dim K_{2n-1,\bV}=\dim \bV^{-1}-\dim \bV^{2n-1}$). Assume $x^{-1}\in \text{Rad}((,)_{-1}|_{K_{2n-1,\bV}})$. Then $\beta(Ax^{-1},K_{2n-1,\bV})=0$ and $\beta(Ax^{-1},I_{2n-1,\bV})=\beta(Ax^{-1},A^{n-1}\bV^{-2n+1})=\beta(A^nx^{-1},\bV^{-2n+1})=0$ (since $x^{-1}\in K_{2n-1,\bV}$). Thus  $\beta(Ax^{-1},\bV^{-1})=0$. It follows that $Ax^{-1}=0$ and  $x^{-1}=0$ (since $A:\bV^{-1}\to \bV^1$ is an isomorphism). Hence (b) holds.

\subsection{}\label{sec-charnot2} Assume in this subsection that $p\neq 2$. Let $\beta_\xi\in\LS(\bV)_2$ and let $A:\bV^a\rightarrow \bV^{a+2} ,\ a\neq -2$ be as in subsection \ref{sec-def}. We define $A:\bV^{-2}\to\bV^0$ by
$\beta(Ax^{-2},x^{0})=\beta_\xi(x^{-2},x^{0}), \forall\ x^0\in \bV^0,x^{-2}\in \bV^{-2}$ (note that $\beta|_{\bV^0}$ is nondegenerate when $p\neq 2$).
The collection of maps $A:\bV^a\rightarrow \bV^{a+2}$ gives rise to a map
$A:\bV\rightarrow\bV,\ A(\sum_ax^a)=\sum_a A x^a.$
Note that for any $v,v'\in\bV$, $\beta(A v,v')=\beta_\xi(v,v').$
It follows that $A\in\Lo(\bV)_2:=\{T\in\Lo(\bV)|T\bV^a\subset\bV^{a+2}\}$.

We show that condition (a) in subsection \ref{sec-def} is equivalent to

(a$'$) $A^{2n}:\bV^{-2n}\to \bV^{2n}$ is an isomorphism for all $n\geq 1$. \\[10pt]
Let $\ker(A^n:\bV^0\rightarrow \bV^{2n})=K_{2n,\bV}$ and $\Ima(A^n:\bV^{-2n}\rightarrow \bV^{0})=I_{2n,\bV}$.
Assume (a$'$) holds. It is clear that all maps $A:\bV^{2n}\to\bV^{2n+2}$ are surjective. We have $\bV^0=K_{2n,\bV}\oplus I_{2n,\bV}$ and thus $\text{Rad}(Q|_{K_{2n,\bV}})\subset\text{Rad}(Q|_{\bV^0})=\{0\}$ (note that $\beta(K_{2n,\bV}, I_{2n,\bV})=0$). Hence (a) holds. Conversely assume (a) holds. Suppose $A^{2n}x^{-2n}=0$. Then $A^nx^{-2n}\in\text{Rad}(Q|_{K_{2n-1,\bV}})=\{0\}$. It follows that  $\beta(x^{-2n},\bV^{2n})=\beta(A^nx^{-2n},\bV^0)=0$ and thus $x^{-2n}=0$.  Hence $A^{2n}:\bV^{-2n}\to\bV^{2n}$ is injective and (a$'$) follows.

It is easy to see that $g\in Z_G(\xi)$ if and only if $g\in Z_G(A)$, $\xi\in\Lg^{*\delta}_2$ if and only if $A\in\Lg^{\delta}_2=\Lo(\bV)_2$, and thus $\xi\in\Lg^{*\delta!}_2$ if and only if $A\in\Lg_2^{\delta!}:=\{T\in\Lg_2^\delta|Z_G(T)\subset G_{\geq 0}^\delta\}$. By \cite[1.5]{Lu1}, $A\in\Lg_2^{\delta!}$ if and only if $A$ satisfies (a$'$) and \ref{sec-def} (b$'$). Proposition \ref{prop-2} follows in this case.

\subsection{}\label{sec-1} We assume $p=2$ in the remainder of this section. Let $\beta_\xi\in\LS(\bV)_2$.
Let $A:\bV^a\rightarrow \bV^{a+2}$, $a\neq -2$ be defined for $\beta_\xi$ as in subsection \ref{sec-def}. Note that $R\subset\bV^0$. Let $\brm$ be the unique integer such that
$$A^{\brm} R\neq0,A^{\brm+1}R=0.$$ We define
$\bv_{\brm}\in R$ by $Q(\bv_{\brm})=1$ and set (if $\brm>0$)
$$\bv_i={A}^{\brm-i}\bv_{\bar{m}},\  i\in[0,\brm-1].$$
Assume $\brm=0$. Let
\begin{eqnarray*}
&& \bbW^0 \text{ be a complementary subspace of }\spn\{\bv_0\}\text{ in }\bV^0\text{ and }\bbW^a=\bV^a\text{ for all }a\neq 0.
 \end{eqnarray*}Assume $\brm>0$. We choose $\bu_0\in \bV^{-2\brm}$ such that $\beta(\bu_0,\bv_0)=1$, and
set
\begin{eqnarray*}
&&\bu_i=A^i\bu_0,\ i\in[1,\brm-1];\\
&&\bbW^{-2i}=\{v\in \bV^{-2i}|\beta(v,\bv_{\brm-i})=0\},\ \bbW^{2i}=\{v\in \bV^{2i}|\beta(v,\bu_{\brm-i})=0\},\ i\in[1,\brm],\\&&\bbW^0=\{v\in \bV^0|\beta_\xi(v,\bu_{\brm-1})=0\},\ \bbW^{a}=\bV^{a}, \text{ if }a\text{ is odd}, \text{ or } a\text{ is even and }a\notin[-2\brm,2\brm].
\end{eqnarray*}
Note that $\bV^{2j}=\bbW^{2j}\oplus\text{span}\{\bv_{\brm-j}\}$, $j\in[0, \brm]$; $\bV^{-2j}=\bbW^{-2j}\oplus\text{span}\{\bu_{\brm-j}\}$, $j\in[1,\brm]$.

We show that $A(\bbW^a)\subset \bbW^{a+2}$, $a\neq -2$. This is clear except when $\brm>0$ and $a=2j$, $j\in[-\brm-1,-2]$ or $j\in[0,\brm-1]$.  If $w^{0}\in \bbW^{0}$, then $\beta(Aw^{0},\bu_{\brm-1})=\beta_\xi(w^{0},
\bu_{\brm-1})=0$; if $w^{2j}\in \bbW^{2j}$, $j\in[1,\brm-1]$, then $\beta(Aw^{2j},\bu_{\brm-j-1})=\beta_\xi(w^{2j},
\bu_{\brm-j-1})=\beta(w^{2j},\bu_{\brm-j})=0$;  if $w^{-2j}\in \bbW^{-2j}$, $j\in[2,\brm]$, then $\beta(Aw^{-2j},\bv_{\brm-j+1})=\beta_\xi(w^{-2j},\bv_{\brm-j+1})=
\beta(w^{-2j},\bv_{\brm-j})=0$; if $w^{-2\brm-2}\in \bbW^{-2\brm-2}$, then $\beta(Aw^{-2\brm-2},\bv_{0})=\beta_\xi(w^{-2\brm-2},\bv_{0})=
\beta(w^{-2\brm-2},A\bv_{0})=0$. Hence in each case $Aw^{2j}\in \bbW^{2j+2}$.

Let $\bar{A}=A|_{\bbW^a}: \bbW^{a}\rightarrow\bbW^{a+2},\ a\neq -2,$
and define $\bar{A}:\bbW^{-2}\rightarrow \bbW^0$ by $\beta(\brA w^{-2},w^0)=\beta_\xi(w^{-2},w^0)$ (note that $\beta|_{\bbW^0}$ is nondegenerate). Then we have a collection of maps $$\bar{A}:\bbW^a\rightarrow \bbW^{a+2}\quad\forall\ a$$
such that $$\beta(\brA w^a,w^{-a-2})=\beta_\xi(w^a,w^{-a-2}), \forall\  w^a\in\bbW^a,\ w^{-a-2}\in\bbW^{-a-2}.$$

\subsection{}\label{sec-equiv}Let $\beta_\xi$,  $\brm$, $\bbW^a$ and $\brA:\bbW^a\to\bbW^{a+2}$ be as in subsection \ref{sec-1}. We show  in this subsection that  condition (a) in subsection \ref{sec-def} holds if and only if the following conditions (a1) and (a2) hold:

(a1) for any $n\in[1,\brm]$, the map $\brA^n:\bbW^{-2n}\to \bbW^0$ is injective and $Q|_{\Ima(\brA^n:\bbW^{-2n}\to \bbW^0)}$ is nondegenerate;\\[5pt]
\indent(a2) for any $n\geq \brm+1$, $\brA^{2n}:\bbW^{-2n}\xrightarrow{} \bbW^{2n}$ is an isomorphism.\\[8pt]
Note that  (a2) implies  that  $\dim \bbW^{2n}$ is even for any $n\geq \brm+1$. In fact it follows from (a2) that $\text{Rad}(Q|_{\text{Im}(\brA^n:\bbW^{-2n}\to \bbW^0)})=\Ima(\brA^{n}:\bbW^{-2n}\to \bbW^{0})\cap\ker(\brA^n:\bbW^0\to \bbW^{2n})=\{0\}$ and thus $\dim \bbW^{-2n}=\dim\text{Im}(\brA^n:\bbW^{-2n}\to \bbW^0)$ is even.

 We first show that \ref{sec-def} (a) holds if and only if (a1$'$) and (a2) hold, where \\[5pt]
\indent(a1$'$) for any $n\in[1,\brm]$,  the map $\bar{A}^n:\bbW^0\xrightarrow{}\bbW^{2n}$ is surjective and $Q|_{\ker(\brA^n:\bbW^0\rightarrow \bbW^{2n})}$ is nondegenerate.

We denote $\ker(A^n:\bV^0\rightarrow \bV^{2n})=K_{2n,\bV},\ \ker(\brA^n:\bbW^0\rightarrow \bbW^{2n})=K_{2n,\bbW}$ and $\Ima(\brA^n:\bbW^{-2n}\rightarrow \bbW^{0})=I_{2n,\bbW}.$
Then one easily shows that
$$K_{2n,\bV}=K_{2n,\bbW},n\in[1,\brm];\ K_{2n,\bV}=K_{2n,\bbW}\oplus\spn\{\bv_{\brm}\}, n\geq \brm+1.$$
\indent Assume \ref{sec-def} (a) holds. We first show that $\brA:\bbW^{2n}\rightarrow \bbW^{2n+2}$ is surjective for all $n\geq 0$. This is clear for $n>\brm$.  Let  $n\in[0,\bar{m}]$ and $w^{2n+2}\in \bbW^{2n+2}$. There exists $v^{2n}\in\bV^{2n}$ such that $Av^{2n}=w^{2n+2}$. We have $\beta(v^{2n},\bu_{\brm-n})=\beta(Av^{2n},\bu_{\brm-n-1})=0$ (if $n\in[0,\brm-1]$) and $Av^{2\brm}=A(w^{2\brm}+a\bv_0)=\brA w^{2\brm}$. Now (a1$'$) follows since $K_{2n,\bV}=K_{2n,\bbW}$ for $n\in[1,\brm]$. We verify (a2).  It is enough to show that $\brA^{2n}:\bbW^{-2n}\xrightarrow{} \bbW^{2n}$ is injective for all  $n\geq \brm+1$. Assume $\brA^{2n}w^{-2n}=0$ for some $w^{-2n}\in\bbW^{-2n}$. Then $\brA^nw^{-2n}\in\text{Rad}(Q|_{K_{2n,\bbW}})=\{0\}$ (note that $\text{Rad}(Q|_{ K_{2n,\bV}})=\text{Rad}(Q|_{ K_{2n,\bbW}})\oplus\spn\{\bv_{\brm}\}$ and that $Q|_{K_{2n,\bV}}$ is nondegenerate). Thus  $w^{-2n}=0$ since $\beta(w^{-2n},\bbW^{2n})=\beta(\brA^nw^{-2n},\bbW^0)=0$ (note that $\beta:\bbW^{-j}\times\bbW^j\to\tk$ is nondegenerate).

Assume (a1$'$) and (a2) hold. It is clear that $A:\bV^{2n}\rightarrow \bV^{2n+2}$ is surjective for all $n\geq 0$, and $Q|_{K_{2n,\bV}}$ is nondegenerate for all $n\in[1,\brm]$ (since $K_{2n,\bV}=K_{2n,\bbW}$). Assume $n\geq \brm+1$. We have $\bbW^0=I_{2n,\bbW}\oplus K_{2n,\bbW}$ (since $I_{2n,\bbW}\cap K_{2n,\bbW}=\{0\}$, $\dim I_{2n,\bbW}=\dim\bbW^{-2n}$ and $\dim K_{2n,\bbW}=\dim\bbW^{0}-\dim \bbW^{2n}$) and thus $\text{Rad}(Q|_{K_{2n,\bbW}})\subset\text{Rad}(Q|_{\bbW^0})=\{0\}$. Hence $Q|_{K_{2n,\bV}}$ is nondegenerate. Thus (a) holds.

It remains to show that (a1) is equivalent to (a1$'$). Note that $\text{Rad}(Q|_{I_{2n,\bbW}})=I_{2n,\bbW}\cap K_{2n,\bbW}$. Assume (a1$'$) holds. Then $\brA^n:\bbW^{-2n}\to \bbW^0$ is injective for any $n\in[1,\brm]$ (if $\brA^nw^{-2n}=0$, then $\beta(w^{-2n},\bbW^{2n})=\beta(w^{-2n},\brA^n\bbW^{0})=\beta(\brA^nw^{-2n},\bbW^{0})=0$ and thus $w^{-2n}=0$). If $I_{2n,\bbW}\cap K_{2n,\bbW}=\{0\}$, then $\bbW^0=I_{2n,\bbW}\oplus K_{2n,\bbW}$ and $\text{Rad}(Q|_{K_{2n,\bbW}})=0$. If $I_{2n,\bbW}\cap K_{2n,\bbW}\neq\{0\}$, then it is contained in $\text{Rad}(Q|_{K_{2n,\bbW}})$. But $\dim\text{Rad}(Q|_{K_{2n,\bbW}})\leq1$. Thus in any case, we have
$\text{Rad}(Q|_{K_{2n,\bbW}})=I_{2n,\bbW}\cap K_{2n,\bbW}=\text{Rad}(Q|_{I_{2n,\bbW}}).$
 Hence (a1) holds.

Conversely assume (a1) holds. Suppose there exists $n\in[1,\brm]$ such that $\brA^n:\bbW^0\to \bbW^{2n}$ is not surjective. We have $\dim K_{2n,\bbW}>\dim \bbW^0-\dim \bbW^{2n}$, $\dim I_{2n,\bbW}=\dim \bbW^{2n}$ and $\dim I_{2n,\bbW}\cap K_{2n,\bbW}=\dim\text{Rad}(Q|_{I_{2n,\bbW}})\leq 1$. Hence $\dim(I_{2n,\bbW}+K_{2n,\bbW})\geq \dim \bbW^0$. It follows that $I_{2n,\bbW}\cap K_{2n,\bbW}\neq\{0\}$ and $\bbW^0=I_{2n,\bbW}+K_{2n,\bbW}$. But then $I_{2n,\bbW}\cap K_{2n,\bbW}\subset\text{Rad}(Q|_{\bbW^0})=\{0\}$ which is a contradiction. Hence $\brA^n:\bbW^0\to \bbW^{2n}$ is  surjective for all $n\in[1,\brm]$. It remains to show that $\text{Rad}(Q|_{K_{2n,\bbW}})=I_{2n,\bbW}\cap K_{2n,\bbW}$.
Let $U^0=I_{2n,\bbW}+K_{2n,\bbW}\subset\bbW^0$. It is easy to see that $\text{Rad}(Q|_{K_{2n,\bbW}})=\text{Rad} (Q|_{U^0})$. We have $\dim(I_{2n,\bbW}\cap K_{2n,\bbW})=\dim \text{Rad}(Q|_{I_{2n,\bbW}})\leq 1$ and thus $\dim U^0\geq\dim \bbW^0-1$. Since $\text{Rad} (Q|_{\bbW^0})=\{0\}$, $\dim\text{Rad} (Q|_{U^0})\leq 1$. Now $I_{2n,\bbW}\cap K_{2n,\bbW}\subset \text{Rad} (Q|_{U^0})$ and thus $I_{2n,\bbW}\cap K_{2n,\bbW}=\text{Rad} (Q|_{U^0})$.  The proof is completed.
\subsection{}\label{sec-cen}
Let $\bbW=\oplus_a\bbW^a$. We have $$\bV=\text{span}\{\bv_i,i\in[0,\brm]\}\oplus\spn\{\bu_i,i\in[0,\brm-1]\}
\oplus\bbW$$ and if $\brm>0$, then \begin{eqnarray*}&&\bbW=\{v\in \bV|\beta(v,\bu_i)=\beta_\xi(v,\bu_i)=0,i\in[0,\brm-1];\ \beta(v,\bv_j)=\beta_\xi(v,\bv_j)=0,j\in[0,\brm]\}.\end{eqnarray*}
The collection of maps $\brA:\bbW^a\rightarrow \bbW^{a+2}$ gives rise to a map $$\brA:\bbW\rightarrow \bbW,\ \brA(\sum_aw^a)=\sum_a \brA w^a.$$
Note that for any $w\in\bbW$ and any $v\in\bV$, we have
$$\beta(\brA w,v)=\beta_\xi(w,v).$$
 It follows that $\brA\in\Lo(\bbW)$.

Let $\pi_{\bbW}:\bV\rightarrow \bbW$ be the natural projection.
Let $g\in Z_G(\xi)=Z_G(\beta_\xi)=\{g\in G|\beta_\xi(gv,gv')=\beta_\xi(v,v')\ \forall\ v,v'\in\bV\}$. We show that
{\em \begin{eqnarray*}
\mathrm{(a)}&&\qquad g\bv_i=\bv_i,\ i\in[0,\brm];\ \pi_{\bbW}(g\brA^jw)=\brA^j\pi_{\bbW}(gw), \text{ and if }\brm>0, \text{ then}\nonumber\\
&&\qquad g\bu_i=\bu_i+\sum_{j\in[0,{\brm}]}b_{j+i}\bv_i+\brA^i\pi_{\bbW}(\bu_0),\ i\in[0,\brm-1];\\&&\qquad g\brA^jw=\sum_{i\in[0,\brm]}\beta(\brA^{i+j}\pi_{\bbW}(\bu_0),\pi_{\bbW}(gw))\bv_i
+\brA^j\pi_{\bbW}(gw),\ w\in\bbW.
\end{eqnarray*}}
We have $g|_R=1$. Thus $g\bv_{\brm}=\bv_{\brm}$. We have $\beta_\xi(g\bv_{\brm},v)=\beta_\xi(\bv_{\brm},g^{-1}v)=\beta(\bv_{\brm-1},g^{-1}v)
=\beta(g\bv_{\brm-1},v)$ and $\beta_\xi(g\bv_{\brm},v)=\beta_\xi(\bv_{\brm},v)=\beta(\bv_{\brm-1},v)$, for any $v\in\bV$. Then $g\bv_{\brm-1}-\bv_{\brm-1}\in R$. Since $Q(g\bv_{\brm-1})=Q(\bv_{\brm-1})$, $g\bv_{\brm-1}=\bv_{\brm-1}$. Similarly, $g\bv_i=\bv_i$, $i\in[0,\brm-2]$.

Now for any $\tilde{w}\in \bbW$, $\beta(\pi_{\bbW}(g\brA w),\tilde{w})=\beta(g\brA w,\tilde{w})=\beta(\brA w,g^{-1}\tilde{w})=\beta_\xi(w,g^{-1}\tilde{w})=
\beta_\xi(gw,\tilde{w})=\beta_\xi(\pi_{\bbW}(gw),\tilde{w})=\beta(\brA\pi_{\bbW}(gw),\tilde{w})$.
Hence
$\pi_{\bbW}(g\brA w)=\brA\pi_{\bbW}(gw)$.

Assume $\brm>0$ and $g\bu_i=\sum_{j\in[0,\brm-1]} a_{ij}\bu_j+\sum_{j\in[0,\brm]}b_{ij}\bv_j+w_i$, where $w_i\in \bbW$. We have $a_{ij}=\beta(g\bu_i,\bv_j)=\beta(g\bu_i,g\bv_j)=\beta(\bu_i,\bv_j)=\delta_{ij}$, $j\in[0,\brm-1]$ and
$\beta_\xi(g\bu_i,v)=\beta_\xi(\bu_i,g^{-1}v)=
\beta(\bu_{i+1},g^{-1}v)=\beta(g\bu_{i+1},v)$. Hence $b_{i,j-1}=b_{i+1,j}$ and $w_{i+1}=\bar{A}w_i$. Assume $w\in \bbW$ and $gw=\sum_{j\in[0,\brm-1]} x_{j}\bu_j+\sum_{j\in[0,\brm]}y_{j}\bv_j+w'$, where $w'\in \bbW$.
Then
$\beta(gw,g\bv_j)=\beta(w,\bv_j)=\beta(gw,\bv_j)=x_j=0$ and thus $\beta(gw,g\bu_i)=y_i+\beta(\brA^iw_0,w')=0$, $i\in[0,\brm-1]$, $\beta_\xi(gw,g\bu_{\brm-1})=y_{\brm}+\beta_\xi(\brA^{\brm-1}w_0,w')=0$.  This completes the proof of (a).

\subsection{} Assume $\beta_\xi\in\LS(\bV)_2^0$ (see \ref{sec-def}). Let $\bar{m}$, $\bbW^a$, $\bbW$ and $\brA$ be defined for $\beta_\xi$ as in subsections \ref{sec-1} and \ref{sec-cen}. Let $\bV^{\geq a}=\oplus_{a'\geq a}\bV^{a'}$ and $\bbW^{\geq a}=\oplus_{a'\geq a}\bbW^{a'}$. We show by induction on $\dim \bV$ that  $Z_G(\xi)\subset G_{\geq 0}^\delta=\{g\in G|g\bV^{\geq a}=\bV^{\geq a}\text{ for all }a\}$ and thus $\xi\in\Lg^{*\delta!}_2$.

 Let $k$ be the largest integer such that $\bV^k\neq 0$. If $k=0$, then $\beta_\xi=0$, $\xi=0$,  $G_{\geq 0}^\delta=G$ and it is clear that $Z_G(\xi)\subset G_{\geq 0}^\delta$. Assume $k\geq 1$. Note that $k\geq 2\brm$ (since $\bv_0\in\bV^{2\brm}$), and for any $w\in \bbW=\bbW^{\geq -k}$, $\brA^{k+1}w\in\bbW^{\geq k+2}=\{0\}$.

Let $g\in Z_G(\xi)$. We first show that $g\bV^k=\bV^k$ and $g\bV^{\geq -k+1}=\bV^{\geq -k+1}$ (note that one follows from the other). Suppose that $k>2\brm$. Then $\bV^k=\bbW^k=\brA^k\bbW$.
Let $w\in \bbW$ and $\pi_{\bbW}(gw)=w'$. Then $g\brA^kw=a\bv_0+\brA^kw'$ (we use \ref{sec-cen} (a) and that $\brA^{k+1}w'=0$).  Since $a^2=Q(g\brA^{k-\brm}w)=Q(\brA^{k-\brm}w)=0$ (we use \ref{sec-cen} (a) and note that $\brA^{k-\brm}\bbW\subset\bbW^{\geq 1}$), $g\brA^kw=\brA^kw'\in\bbW^k$. It follows that $g\bV^k=\bV^k$. Suppose now that $k=2\brm$. Note that $\brA\in\Lo(\bbW)_2^0$ (with respect to the $o$-good grading $\bbW=\oplus\bbW^a$ of $\bbW$, see \cite[1.5]{Lu1} for the definition of $\Lo(\bbW)_2^0$ and see \ref{sec-def} (b$'$), \ref{sec-equiv} (a1), (a2)). Thus  by \cite[1.8]{Lu1}, $\bbW^{\geq -2\brm+1}=\{x\in \bbW|\brA^{2\brm}x=0,Q(\brA^{\brm}x)=0\}$. We have $\bV^{\geq -2\brm+1}=\spn\{\bv_i,i\in[0,\brm]\}\oplus\spn\{\bu_i,i\in[1,\brm-1]\}\oplus \bbW^{\geq -2\brm+1}$.
 It is clear that $g\bv_i\in\bV^{\geq -2\brm+1}$ (see \ref{sec-cen} (a)). Note that for any $w\in \bbW$ and any $i\geq \brm+1$, $\brA^i w\in\bbW^{\geq 2}$ and thus $Q(\brA^iw)=0$. Hence $ \brA\bbW\subset \bbW^{\geq -2\brm+1}$. It follows that $g\bu_j\in \bV^{\geq -2\brm+1}$, $j\in[1,\bar{m}-1]$ (see \ref{sec-cen} (a)). For any $w\in \bbW^{\geq -2\brm+1}$,  $\brA^{2\brm}\pi_{\bbW}(gw)=\pi_{\bbW}(g\brA^{2\brm}w)=0$ (note that $\brA^{2\brm}w\in\bbW^{\geq 2\bar{m}+1}=\{0\}$),  $Q(\brA^{m}\pi_{\bbW}(gw))
=Q(g\brA^{\brm}w)=Q(\brA^{\brm}w)=0$ (note that $\brA^{\brm}w\in\bbW^{\geq 1}$) and thus $gw\in\bV^{\geq -2\brm+1}$. Hence  $g\bV^{\geq -2\brm+1}=\bV^{\geq -2\brm+1}$.

Let $\bV'=\bV^{\geq-k+1}/\bV^k$. Then $Q$ induces a nondegenerate quadratic form $Q'$ on $\bV'$.  There is a natural ($o$-good) grading $\bV'=\bV^{-k+1}\oplus\cdots\oplus \bV^{k-1}$ on $\bV'$ and $\beta_\xi$ induces a symplectic form $\beta_\xi'\in\LS(\bV')_2^0$ with respect to this grading (note that $\beta_\xi(\bV^k,\bV^{\geq -k+1})=\beta_\xi(\bbW^k,\bbW^{\geq -k+1})=\beta(\brA\bbW^k,\bbW^{\geq -k+1})=0$). Now $g\in Z_G(\xi)$ induces an element $g'\in Z_{G'}(\beta_\xi')$, where $G'=SO(\bV')$ is  defined with respect to $Q'$. By induction hypothesis, the subspace $\bV^{a}\oplus\cdots\bV^{a+1}\oplus\cdots\oplus\bV^{k-1}$ of $\bV'$ is $g'$-stable for all $a\in[-k+1,k-1]$. It follows that $g \bV^{\geq a}=\bV^{\geq a}$ for all $a$ and thus $g\in G_{\geq 0}^\delta$.

\subsection{}\label{sec-2} Assume $\xi\in\Lg^{*\delta!}_2$. We show that $\beta_\xi\in\LS(\bV)_2^0$. Let $\bar{m}$, $\bbW^a$, $\bbW$ and $\brA$ be defined for $\beta_\xi$ as in subsections \ref{sec-1} and \ref{sec-cen}.

We first show that $\dim \bbW^{-2\brm}\geq\dim\bbW^{-2\brm-2}$. Otherwise
 $\brA:\bbW^{-2\brm-2}\to\bbW^{-2\brm}$ is not injective. We choose a nonzero $e^{-2\brm-2}\in \bbW^{-2\brm-2}$ such that $\brA e^{-2\brm-2}=0$. Assume $\brm=0$. Define $g:\bV\to\bV$, $g\notin G_{\geq 0}^\delta$ by
\begin{equation*}
g\bv_{\brm}=\bv_{\brm};\ gw=w+\beta(w,e^{-2\brm-2})(e^{-2\brm-2}+\bv_{\brm}),\   w\in \bbW.
\end{equation*}
Let $v,v'\in\bV$ and $w=\pi_{\bbW}(v)$, $w'=\pi_{\bbW}(v')$. Since $\beta_\xi(\bv_{\brm},\bV)=0$ and $\beta_\xi(e^{-2\brm-2},\bV)=\beta(\brA e^{-2\brm-2},\bV)=0$, it is easy to verify that $Q(gv)=Q(w)+a^2=Q(v)$ and $\beta_\xi(gv,gv')
=\beta_\xi(w,w')
=\beta_\xi(v,v').$
Thus $g\in Z_G(\xi)\nsubseteq G_{\geq 0}^\delta$ which is a contradiction. Assume $\brm>0$. Define $g:\bV\to\bV$, $g\notin G_{\geq 0}^\delta$ by
\begin{eqnarray*}
&&g\bv_i=\bv_i,\ i\in[0,\bar{m}];\ g\bu_0=\bu_0+e^{-2\brm-2}; \ g\bu_i=\bu_i,\ i\in[1,\brm-1];\\
&&gw=w+\beta(w,e^{-2\brm-2})\bv_0,\  w\in \bbW.
\end{eqnarray*}
For any $v=w+\sum  a_i\bu_i+\sum  b_i\bv_i$ and $v'=w'+\sum a_i'\bu_i+\sum b_i'\bv_i$ in $\bV$ (where $w=\pi_{\bbW}(v)$, $w'=\pi_{\bbW}(v')$),  since $\beta_\xi(\bv_0,\bV)=0$ and $\beta_\xi(e^{-2\brm-2},\bV)=\beta(\brA e^{-2\brm-2},\bV)=0$, we have
\begin{eqnarray*}
&Q(gv)=Q(v+\beta(w,e^{-2\brm-2})\bv_0+a_0e^{-2\brm-2})
=Q(v)+\beta(v,\beta(w,e^{-2\brm-2})\bv_0+a_0e^{-2\brm-2})=Q(v),\\&
\beta_\xi(gv,gv')=\beta_\xi(v+\beta(w,e^{-2\brm-2})\bv_0+a_0e^{-2\brm-2},
v'+\beta(w',e^{-2\brm-2})\bv_0+a_0'e^{-2\brm-2})=\beta_\xi(v,v').\end{eqnarray*}  Thus $g\in Z_G(\xi)\nsubseteq G_{\geq 0}^\delta$, which is again a contradiction.

Hence $\bbW=\oplus_a\bbW^a$ is an $o$-good grading of $\bbW$. Note that $\brA\in\Lo(\bbW)_2=\{T\in\Lo(\bbW)|T\bbW^a\subset\bbW^{a+2}\}$. Let $SO(\bbW)_{\geq 0}=\{g\in SO(\bbW)|g\bbW^{\geq a}=\bbW^{\geq a}, \ \forall\ a \}$.

 We show that $\brA\in \Lo(\bbW)_2^0$,  namely, $\brA^i:\bbW^{-i}\rightarrow \bbW^i$ is an isomorphism for any odd $i$,  $\brA^{i}:\bbW^{-2i}\to \bbW^0$ is injective and $Q|_{\text{Im}(\brA^i:\bbW^{-2i}\to \bbW^0)}$ is nondegenerate for all $i\geq 1$ (see \cite[1.5]{Lu1}).  Otherwise by \cite[1.8]{Lu1} there exists $g_0\in Z_{SO(\bbW)}(\brA)$ such that $g_0\notin SO(\bbW)_{\geq 0}$. Define $g:\bV\to \bV$ by $gw=g_0w$ for all $w\in \bbW$ and $g\bv_i=\bv_i$, $g\bu_i=\bu_i$. Note that $g\notin G_{\geq 0}^\delta$. For any $v,v'\in\bV$ with $\pi_{\bbW}(v)=w$ and $\pi_{\bbW}(v')=w'$, since $\beta_\xi(g_0w,g_0w')=\beta(\brA g_0w,g_0w')=\beta(g_0\brA w,g_0w')=\beta(\brA w,w')=\beta_\xi(w,w')$, we have
\begin{eqnarray*}&&Q(gv)=Q(v+g_0w+w)=Q(v)+\beta(v,g_0w+w)+Q(g_0w+w)\\&&\qquad
=Q(v)+\beta(w,g_0w+w)+Q(g_0w)+Q(w)+\beta(g_0w,w)=Q(v)\\
&&\beta_\xi(gv,gv')=\beta_\xi(v+g_0w+w,v'+g_0w'+w')=\beta_\xi(v,v')+\beta_\xi(w,g_0w'+w')+\beta_\xi(g_0w+w,w')
\\
&&\quad+\beta_\xi(g_0w+w,g_0w'+w')=\beta_\xi(v,v')+\beta_\xi(g_0w,g_0w')+\beta_\xi(w,w')
=\beta_\xi(v,v').
\end{eqnarray*}
Thus $g\in Z_G(\xi)\nsubseteq G_{\geq 0}^\delta$ which is a contradiction.

We show that for all $k\geq \brm+1$, $\brA^{2k}:\bbW^{-2k}\rightarrow \bbW^{2k}$ is an isomorphism. Otherwise there exist a $k$ and a nonzero $e^{-2k}\in \bbW^{-2k}$ such that $\brA^{2k}e^{-2k}=0$. We have $e^0=\brA^ke^{-2k}\neq 0$ (since $\brA^k:\bbW^{-2k}\to\bbW^0$ is injective). Let $e^{2j-2k}=\brA^je^{-2k}$ for $j\in[0,2k]$. Note that  $\beta(e^{2j-2k},e^{2k-2j})=\beta(\brA^{2k}e^{-2k},e^{-2k})=0$ and $\beta_\xi(e^{2j-2k},e^{2k-2j-2})=\beta(\brA e^{2j-2k},e^{2k-2j-2})=\beta(e^{2j-2k+2},e^{2k-2j-2})=0$. Fix a square root $\sqrt{Q(e^0)}$ of $Q(e^0)$. Define $g:\bV\to\bV$, $g\notin G_{\geq 0}^\delta$ by
\begin{eqnarray*}
&&g\bv_i=\bv_i,\ i\in[0,\brm];\ g \bu_i=\bu_i+\sqrt{Q(e^0)}e^{2i-2\brm-2},\ i\in[0,\brm-1];\\
&&g w=w+\sum_{j\in[0,2k-1]}\beta(e^{2k-2j-2},w)e^{2j-2k}+
\sum_{j\in[0,\brm]}\beta(e^{2j-2\brm-2},w)\sqrt{Q(e^0)}\bv_{j},\ w\in\bbW.
\end{eqnarray*}
For any $v=w+\sum a_i\bv_i+\sum b_i\bu_i$ and $v'=w'+\sum a_i'\bv_i+\sum b_i'\bu_i$ in $\bV$ (where $w=\pi_{\bbW}(v)$, $w'=\pi_{\bbW}(v')$), we have
\begin{eqnarray*}
&&Q(gv)=Q(v+\sum_{i\in[0,\brm-1]}b_i\sqrt{Q(e^0)}e^{2i-2\brm-2}+\sum_{j\in[0,2k-1]}\beta(e^{2k-2j-2},w)e^{2j-2k}\\&&\qquad
+
\sum_{j\in[0,\brm]}\beta(e^{2j-2\brm-2},w)\sqrt{Q(e^0)}\bv_j)\\&=&Q(v)+\sum_{i\in[0,\brm-1]}\sqrt{Q(e^0)}\beta(w,e^{2i-2\brm-2})b_i+\sum_{j\in[0,2k-1]}\beta(e^{2k-2j-2},w)\beta(w,e^{2j-2k})
\\&&\qquad+
\sum_{j\in[0,\brm-1]}\beta(e^{2j-2\brm-2},w)\sqrt{Q(e^0)}b_j+\beta(e^{-2},w)^2Q(e^0)+(\beta(e^{-2},w)
\sqrt{Q(e^0)})^2\\
&=&Q(v)+\sum_{j\in[0,2k-1]}\beta(e^{2k-2j-2},w)\beta(w,e^{2j-2k})=Q(v)
\\&&\beta_\xi(gv,gv')=\beta_\xi(v+\sum_{i\in[0,\brm-1]}b_i\sqrt{Q(e^0)}e^{2i-2\brm-2}
+\sum_{j\in[0,2k-1]}\beta(e^{2k-2j-2},w)e^{2j-2k}\\&&\qquad+
\sum_{j\in[0,\brm]}\beta(e^{2j-2\brm-2},w)\sqrt{Q(e^0)}\bv_j,v'+
\sum_{i\in[0,\brm-1]}b_i'\sqrt{Q(e^0)}e^{2i-2\brm-2}\\&&\qquad
+\sum_{j\in[0,2k-1]}\beta(e^{2k-2j-2},w')e^{2j-2k}+
\sum_{j\in[0,\brm]}\beta(e^{2j-2\brm-2},w')\sqrt{Q(e^0)}\bv_j)\end{eqnarray*}
\begin{eqnarray*}
&=&\beta_\xi(v,v')+\sum_{i\in[0,\brm-1]}b_i\sqrt{Q(e^0)}\beta_\xi(e^{2i-2\brm-2},w')
+\sum_{j\in[0,2k-1]}\beta(e^{2k-2j-2},w)\beta_\xi(w',e^{2j-2k})\\&&
\qquad+\sum_{j\in[1,\brm]}\beta(e^{2j-2\brm-2},w)\sqrt{Q(e^0)}b_{j-1}'+
\sum_{i\in[0,\brm-1]}b_i'\sqrt{Q(e^0)}\beta_\xi(e^{2i-2\brm-2},w)\\&&\qquad
+\sum_{j\in[0,2k-1]}\beta(e^{2k-2j-2},w')\beta_\xi(e^{2j-2k},w)+
\sum_{{j\in[1,\brm]}}\beta(e^{2j-2\brm-2},w')\sqrt{Q(e^0)}b_{j-1}\end{eqnarray*}
\begin{eqnarray*}
&=&\beta_\xi(v,v')+\sum_{i\in[0,\brm-1]}b_i\sqrt{Q(e^0)}\beta(e^{2i-2\brm},w')
+\sum_{j\in[0,2k-2]}\beta(e^{2k-2j-2},w)\beta(w',e^{2j-2k+2})\\
&&
\qquad+\sum_{{j\in[1,\brm]}}\beta(e^{2j-2\brm-2},w)\sqrt{Q(e^0)}b_{j-1}'+
\sum_{i\in[0,\brm-1]}b_i'\sqrt{Q(e^0)}\beta(e^{2i-2\brm},w)\\&&\qquad+\sum_{j\in[0,2k-2]}\beta(e^{2k-2j-2},w')\beta(e^{2j-2k+2},w)+
\sum_{{j\in[1,\brm]}}\beta(e^{2j-2\brm-2},w')\sqrt{Q(e^0)}b_{j-1}\\
&=&\beta_\xi(v,v').
\end{eqnarray*}
Thus $g\in Z_G(\xi)\nsubseteq G_{\geq 0}^\delta$ which is again a contradiction.  This completes the proof of Proposition \ref{prop-2}.

\section{Proof of Theorem \ref{thm}}\label{sec-thm}
\subsection{}\label{ssec-ind}

For $\bV_*=(\bV^{\geq a})\in\mathfrak{F}_o(\bV)$ (see \ref{sec-ogood}), let $\eta(\bV_*)$ be the set of all $\beta_\xi\in\mathfrak{S}(\bV)$ such that $\beta_\xi(\bV^{\geq a},\bV^{\geq b})=0$ whenever $a+b\geq -1$, and that the symplectic bilinear form $\bbeta_\xi$ induced by $\beta_\xi$ is in $\LS(\bV)_{ 2}^0$ (with respect to the corresponding $o$-good grading $\bV=\oplus_{a\in\mathbb{Z}}\bV^a$ such that $\bV^{\geq a}=\oplus_{a'\geq a}\bV^{a'}$), where ${\bbeta}_\xi\in\LS(\bV)_2$ is defined as follows
\begin{equation*}
{\bbeta}_\xi(\sum_a x^a,\sum_a y^a)=\sum_a\beta_\xi(x^a,y^{-a-2}), \text{ for all }x^a\in \bV^a,y^a\in\bV^a.
\end{equation*}
Note that $\eta(\bV_*)\subset\LS(\bV)_{nil}.$

\begin{proposition}\label{prop-1}
The map $$\sqcup_{\bV_*\in\mathfrak{F}_o(\bV)}\eta(\bV_*)\rightarrow\LS(\bV)_{nil},\ \beta_\xi\mapsto\beta_\xi$$ is a bijection.
\end{proposition}

When $p\neq 2$, the map $\beta_\xi\mapsto A$, where $\beta_\xi(v,v')=\beta(Av,v')$ for all $v,v'\in\bV$, defines a bijection $\LS(\bV)_{nil}\xrightarrow{\sim}\cN_{\Lg}$  and a bijection $\eta(\bV_*)\xrightarrow{\sim}\eta'(\bV_*)$, where $\cN_\Lg$ is the set of nilpotent elements in $\Lg$ and $\eta'(\bV_*)$ is as in \cite[A.4]{Lu1}. Thus the proposition follows from \cite[A.4(a)]{Lu1} in this case. The proof when $p=2$ will be given in subsections \ref{sec-vu}-\ref{sec-surj}.

In view of the identification of $\mathfrak{F}_o(\bV)$ with $D_G$ (see \ref{sec-ogood}) and the identification of $\LS(\bV)_{nil}$ with $\cN_{\Lg^*}$ (see \ref{ssec-map}), Theorem \ref{thm} follows from Proposition \ref{prop-1} and Proposition \ref{prop-2}.

\subsection{}\label{sec-vu}We assume $p=2$ through subsection \ref{sec-surj}. For any nonzero element $\beta_\xi\in\LS(\bV)_{nil}$ we associate to $\beta_\xi$ integers $m,\lambda_1,l_1$ in $\mathbb{N}$ and a subspace $H_{\beta_\xi}$ of $\bV$ as follows.

 Let $m\in\mathbb{N}$ be the unique integer and $\{v_i,i\in[0,m]\}$ the unique set of vectors defined for $\beta_\xi$ as in \ref{ssec-map} (a).
Assume $m>0$. We choose a vector $u_0\in \bV$ such that

\qquad (a) {\em $\beta(u_0,v_0)=1,\beta(u_0,v_i)=0,i\in[1,m-1],\ Q(u_0)=0.$}\\[10pt]
Let $\{u_i, i\in[1,m-1]\}$ be the unique set of vectors (see \cite[Lemma 3.6]{X2}) such that
 $$Q(u_i)=0;\ \beta(u_i,v)=\beta_\xi(u_{i-1},v),\ \forall\ v\in\bV.$$ Let $\bW$ be a subspace  of $\bV$ such that
 \begin{eqnarray*}
\mathrm{(b)}&&\bW\oplus\spn\{v_m\}=\bV,\text{ if }m=0;\nonumber\\
&&\bW=\{v\in \bV|\beta(v,v_i)=\beta(v,u_i)=0,i\in[0,m-1],\ \beta_\xi(v,u_{m-1})=0\},\text{ if }m>0.
 \end{eqnarray*}
Then $\bV=\spn\{v_i,i\in[0,m]\}\oplus\spn\{u_i,i\in[0,m-1]\}\oplus\bW$ (see \cite[Lemma 3.8]{X2}) and $\beta|_{\bW}$ is nondegenerate. Define
$$T_\xi:\bW\to \bW \text{ by }\beta(T_\xi w,w')=\beta_\xi(w,w'),\ \forall\ w,w'\in \bW.$$
Then for any $x\in \bW$ and any $v\in \bV$, $\beta_\xi(x,v)=\beta(T_\xi x,v)$. Moreover $T_\xi\in\Lo(\bW)$ is nilpotent (see \cite[Lemma 3.11]{X2}). Let $\pi_\bW:\bV\to\bW$ denote the natural projection.

 Let $\lambda_1$ be the smallest integer such that $T_\xi^{\lambda_1}\bW=0$. Let $$l_1=\max(\lambda_1-m,f),$$ where $f$ is the smallest integer such that $Q(T_\xi^{f}\bW)=0$. Let $$\rho:\ker T_\xi^{\lambda_1-1}\to\tk$$ be the map  $w\mapsto Q(T_\xi^{l_1-1}w)$. We set
\begin{eqnarray*}
&&\Lambda_\bW=\{x\in\bW|Q(T_\xi^{l_1-1}x)=0\}\ \text{ if }\lambda_1-l_1<m< l_1,\\
&&\Lambda_\bW=\ker\rho\quad \text{if }m=\lambda_1-l_1<l_1;
\end{eqnarray*}
and
\begin{eqnarray*}
&&H_{\beta_\xi}=\spn\{v_m\}\oplus\Lambda_\bW,\text{ if }m=0;\\
&& H_{\beta_\xi}=\spn\{v_i,i\in[0,m]\}\oplus\spn\{u_i,i\in[1,m-1]\}\oplus \bW,\ \text{ if } m\geq l_1;\\
&&H_{\beta_\xi}=\spn\{v_i,i\in[0,m]\}\oplus\spn\{u_i,i\in[1,m-1]\}\oplus\Lambda_{\bW},\ \text{  if }\lambda_1-l_1<m<l_1,\\
&&\qquad\qquad0<m=\lambda_1-l_1=l_1-1, \text{ or }0<m=\lambda_1-l_1<l_1-1\text{ and }\rho\neq 0;\\
&&H_{\beta_\xi}=\spn\{v_i,i\in[0,m]\}\oplus\spn\{u_i,i\in[1,m-1]\}\oplus \spn\{u_0+w_{**}\}\oplus\Lambda_\bW,\\
&&\qquad\qquad \text{  if }0<m=\lambda_1-l_1<l_1-1\text{ and }\rho=0,
\end{eqnarray*}
where $w_{**}\in\bW$ is defined as follows. There exists a unique $w_*\in \bW$ such that $\beta({w_*,w})^2=Q(T_\xi^{l_1-1}w)$ for all $w\in \bW$. Since $\rho=0$, $w_*\in(\ker T_\xi^{\lambda_1-1})^{\p}\cap\bW=T_\xi^{\lambda_1-1}\bW$.  We  choose any $w_{**}\in \bW$ such that $w_*=T_\xi^{\lambda_1-1}w_{**}$ (the defintion of $H_{\beta_\xi}$ does not depend on the choice). We show that

(c) {\em the definitions of $\lambda_1$, $l_1$ and $H_{\beta_\xi}$ do not depend on the choice of $\bW$ (if $m=0$) and $u_0$ (if $m>0$).}

Assume $m=0$. Let $\tilde{\bW}$ be another complementary subspace of $\spn\{v_m\}$ in $\bV$ and let $\tilde{T_\xi}:\tilde{\bW}\to\tilde{\bW}$ be defined as $T_\xi:\bW\rightarrow \bW$. Then one easily shows that for any $\tilde{w}\in\tilde{\bW}$, $\pi_\bW(\tilde{T}_\xi\tilde{w})=T_\xi(\pi_\bW(\tilde{w}))$. Now (c) follows (note that $\rho=0$ in this case and thus $\Lambda_\bW=\ker T_\xi^{\lambda_1-1}$).

Assume $m>0$. We choose another $\tilde{u}_0$ satisfying (a) and define $\tilde{u}_i$, $i\in[1,m-1]$, $\tilde{\bW}$, $\tilde{T}_\xi:\tilde{\bW}\to\tilde{\bW}$ as $u_i$, $i\in[1,m-1]$, $\bW$, $T_\xi:\bW\to\bW$. Let $\tilde{\lambda}_1,\tilde{l}_1$ be defined for $\tilde{T}_\xi$ and $\tilde{\bW}$ as $\lambda_1,l_1$. By the same argument as in the proof of \ref{sec-cen} (a) one shows that\\[5pt]
\indent\quad (d) {\em $\tilde{u}_j=u_j+\sum_{i\in[0,m]} a_{i+j}v_i+T_\xi^jw_0,\ j\in[0,m-1],$}\\[5pt]
\indent\quad (e) {\em $\tilde{\bW}=\{\sum_{i\in[0,m]}\beta(w,T_\xi^iw_0)v_i+w|w\in \bW\},\ \pi_{\bW}(\tilde{T}_\xi\tilde{w})=T_\xi (\pi_{\bW}(\tilde{w})),\ \tilde{w}\in\tilde{\bW}$,}\\[6pt]
where $w_0=\pi_\bW(\tilde{u}_0)$ and $a_{2j}+a_{m+j}^2+Q(T_\xi^jw_0)=0$, $j\in[0,m-1]$.
 It follows from (e) that $T_\xi^e\bW=0$ if and only if $\tilde{T}_\xi^e\tilde{\bW}=0$ and
 thus $\tilde{\lambda}_1=\lambda_1$. It also follows that $Q(\tilde{T}_\xi^{l_1}\tilde{w})=Q(T_\xi^{l_1}(\pi_\bW(\tilde{w})))=0$ for all $\tilde{w}\in\tilde{\bW}$ (note that $m+l_1\geq\lambda_1$). Thus $\tilde{l}_1\leq l_1$. Similarly $l_1\leq\tilde{l}_1$. Hence $l_1=\tilde{l}_1$.

Let  $\tilde{\rho}:\ker \tilde{T}_\xi^{\lambda_1-1}\to\tk$ be the map  $\tilde{w}\mapsto Q(\tilde{T}_\xi^{l_1-1}\tilde{w})$.  For any $\tilde{w}\in\tilde{\bW}$, $\tilde{T}_\xi^{\lambda_1-1}\tilde{w}=0$ if and only if $T_\xi^{\lambda_1-1}(\pi_\bW(\tilde{w}))=0$; and  if $\tilde{T}_\xi^{\lambda_1-1}\tilde{w}=0$, then
$Q(\tilde{T}_\xi^{l_1-1}\tilde{w})=Q(T_\xi^{l_1-1}(\pi_\bW(\tilde{w})))$ (we use (e)). Thus $\rho=0$ if and only if $\tilde{\rho}=0$.
It also follows that

\quad (f) {\em $\ker\tilde{\rho}\subset\ker\rho\oplus\spn\{v_i,i\in[0,m]\}$.}\\[5pt]
Now we denote $H_{u_0}\subset\bV$ the subspace in the r.h.s. of the definition of $H_{\beta_\xi}$ and let $H_{\tilde{u}_0}\subset\bV$ be defined as $H_{u_0}$ with $u_i,\bW,T_\xi$, $\Lambda_\bW$ replaced by $\tilde{u}_i,\tilde{\bW},\tilde{T}_\xi$, $\Lambda_{\tilde{\bW}}$. We need to show that $H_{u_0}=H_{\tilde{u}_0}$. If $m\geq l_1$, then this  is clear since $H_{\tilde{u}_0}\subset H_{u_0}$ (see (d) and (e)) and $\dim H_{\tilde{u}_0}=\dim H_{u_0}$. Assume $m<l_1$. We show that\\[5pt]
\indent\quad (g) {\em $\pi_\bW(\tilde{u}_j)\in \Lambda_{\bW}$ for all $j\geq 1$;}\\[5pt]
\indent\quad (h) {$\pi_{\bW}(\Lambda_{\tilde{\bW}})\subset\Lambda_{\bW}$.}\\[5pt]
In fact, for all $j\geq 1$,  we have $T_\xi^{\lambda_1-1}(T_\xi^jw_0)=0$,  $Q(T_\xi^{l_1-1}(T_\xi^jw_0))=Q(T^{l_1}(T_\xi^{j-1}w_0))=0$ and thus (g) follows from (d). Now   if $m> \lambda_1-l_1$,
then $Q(\tilde{T}_\xi^{l_1-1}\tilde{w})=Q(T_\xi^{l_1-1}(\pi_\bW(\tilde{w}))$ for any $\tilde{w}\in\tilde{\bW}$ and thus (h) follows; if $m=\lambda_1-l_1$, then (h) follows from (f).

Assume $\lambda_1-l_1<m<l_1$, $m=\lambda_1-l_1=l_1-1$, or $m=\lambda_1-l_1<l_1-1$ and $\rho\neq 0$.  It follows from (g) and (h) that $H_{\tilde{u}_0}\subset H_{u_0}$ (see (e)). Similarly $H_{u_0}\subset H_{\tilde{u}_0}$.

Assume $m=\lambda_1-l_1<l_1-1$ and $\rho=0$.   Suppose that $\tilde{w}_*\in\tilde{\bW}$ satisfies $\beta(\tilde{w}_*,\tilde{w})^2=Q(\tilde{T}_\xi^{l_1-1}\tilde{w})$ for any $\tilde{w}\in\tilde{\bW}$. Then  $\beta(\pi_\bW(\tilde{w}_*),w)^2=Q(T_\xi^{l_1-1}w)+\beta(w,T_\xi^{\lambda_1-1}w_0)^2
=\beta(w_*,w)^2+\beta(w,T_\xi^{\lambda_1-1}w_0)^2$ for all $w\in\bW$ and thus  $\beta(\pi_\bW(\tilde{w}_*)+w_*+T_\xi^{\lambda_1-1}w_0,\bW)=0$. Hence $\pi_\bW(\tilde{w}_*)=w_*+T_\xi^{\lambda_1-1}w_0=T_\xi^{\lambda_1-1}(w_{**}+w_0)$ and we can choose $\tilde{w}_{**}=w_{**}+w_0+\sum \beta(w_{**},T_\xi^iw_0)v_i$. It follows
that $\tilde{u}_0+\tilde{w}_{**}\in H_{u_0}$. It then follows from (g) and (h) that ${H}_{\tilde{u}_0}\subset H_{u_0}$. Similarly ${H}_{{u_0}}\subset H_{\tilde{u}_0}$.
This completes the proof of (c).

Let $L=H_{\beta_\xi}^\p\cap Q^{-1}(0)$. We show that

\quad (i) {\em $\beta_\xi(L,\bV)=0.$}\\[10pt]
If $m\geq l_1$, then $m>0$ and $L=\spn\{v_0\}$.
If $\lambda_1-l_1<m< l_1$, $m=\lambda_1-l_1<l_1-1$ and $\rho\neq 0$, or $m=\lambda_1-l_1=l_1-1$, then $m>0$ and $L=\spn\{v_0\}\oplus (\Lambda_{\bW}^{\p}\cap\bW)$.
If $m=\lambda_1-l_1<l_1-1$ and $\rho=0$, then either $m=0$ or $m>0$ in which case we can choose $u_0$ such that $w_{**}=0$ (see the definition of $H_{\beta_\xi}$), and thus $L= \Lambda_{\bW}^{\p}\cap\bW.$ It is easy to see that  $T_\xi \bW\subset\Lambda_\bW$. Hence for any $x\in\Lambda_{\bW}^{\p}\cap\bW$, $\beta_\xi(x,\bW)=\beta(T_\xi x,\bW)=\beta(x,T_\xi \bW)=0$. Now (i) follows since $\beta_\xi(v_0,\bV)=0$ when $m>0$.

\subsection{} \label{ssec-vn}Let $\bV=\oplus_a\bV^a$ be an $o$-good grading such that $\beta_\xi\in \eta(\bV_*)$, where $\bV_*=(\bV^{\geq a})$ and $\bV^{\geq a}=\oplus_{a'\geq a}\bV^{a'}$. Let
$\bar{\beta}_\xi\in\LS(\bV)_2^0$ be the symplectic bilinear form induced by $\beta_\xi$ (see \ref{ssec-ind}). Let  $H_{\beta_\xi}$,  $m$, $\{v_i,i\in[0,m]\}$, $\lambda_1$, $l_1$ be defined for $\beta_\xi$ as in subsection \ref{sec-vu}. Let $n$ be the largest integer such that $\bV^n\neq 0$. We show in this subsection that

\begin{lem} We have $\bV^{\geq -n+1}=H_{\beta_\xi}$ and $n=\max(2m,m+l_1-1)$.
\end{lem}
In general, for $x\in\bV=\oplus\bV^a$, $x^a$ denotes the $\bV^a$ component of $x$. We have\\[5pt]
\indent\quad (a) {\em $v_m\in \bV^0
\text{ and }\text{if }m>0,\ v_{i}\in \bV^{\geq j},\text{ then }v_{i-1}\in\bV^{\geq j+2},i\in[1,m].$}\\[5pt]
In fact, if $v_{i}\in \bV^{\geq j}$, then $\beta(v_{i-1},\bV^{\geq -1-j})=\beta_\xi(v_{i},\bV^{\geq -1-j})=0$. Hence $v_{i-1}\in (\bV^{\geq -1-j})^\p\cap Q^{-1}(0)=\bV^{\geq j+2}$. It follows that $v_0\in \bV^{\geq 2m}$ and  thus
$$
n\geq 2m.
$$
Let $n_0\in[0, n-2m]$ be the unique integer such that
\begin{equation*}\label{v0}v_0\in \bV^{\geq n-n_0},\ v_0\notin \bV^{\geq n-n_0+1}.\end{equation*} If $m>0$, we (can) choose $u_0$ such that $u_0\in\bV^{\geq -n+n_0}$ (see \ref{sec-vu} (a)). Then the same argument as in the proof of (a) shows that

\quad (b) {\em $u_i\in \bV^{\geq -n+n_0+2i},\ i\in[0,m-1].$}\\[10pt]
It follows from (a) and (b) that

\quad (c) {\em $v_i\in \bV^{\geq -n+1}, i\in[0,m];\ u_i\in \bV^{\geq
-n+1}, i\in[1,m-1];\   \text{ and if }n_0>0,\text{ then }u_0\in\bV^{\geq -n+1}.$}\\[10pt]
Let $\bW\subset\bV$ be chosen as in \ref{sec-vu} (b). We have

\quad (d) {\em
$\text{if }x\in \bW\cap \bV^{\geq a},\text{  then } T_\xi^bx\in\bV^{\geq a+2b}+\text{span}\{v_i,i\in[\max(0,{m-b+1}),m]\}.$}\\[10pt]
In fact, suppose that $T_\xi^{b-1}x=y+\sum_{k\in[\max(0,{m-b+1}),m]}c_kv_{k}$, where $y\in\bV^{\geq a+2b-2}$. Then $\beta(T_\xi^bx,v)=\beta_\xi(T_\xi^{b-1}x,v)=\beta_\xi(\sum_{k\in[\max(0,{m-b+1}),m]}c_kv_{k},v)=\beta(\sum_{k\in[\max(0,{m-b}),m-1]}c_{k+1}v_{k},v)$  $\forall\ v\in\bV^{\geq 1-a-2b}$. Hence $T_\xi^bx-\sum_{k\in[\max(0,{m-b}),m-1]}c_{k+1}v_{k}\in(\bV^{\geq 1-a-2b})^\p=\bV^{\geq a+2b}+R$ and (d) follows.

It follows from (d) that for all $x\in \bW=\bW\cap\bV^{\geq -n}$,  $T_\xi^{n+1}x\in \bW\cap(\bV^{\geq n+2}+\spn\{v_i\})=\{0\}$, and $T_\xi^{\theta}x\in \bW\cap(\bV^{\geq -n+2\theta}+ \text{span}\{v_i\})$, where $\theta=[\frac{2n-2m-n_0+2}{2}]$. Assume $T_\xi^{\theta}x=x_1+\sum_{i=0}^m c_i v_{i}$, where $x_1\in\bV^{\geq -n+2\theta}$.
Then $c_m=\beta_\xi(T_\xi^{\theta}x,u_{m-1})=0$ (note that  $u_{m-1}\in \bV^{\geq -n+n_0+2m-2}$ and thus $\beta_\xi(x_1,u_{m-1})=0$). It follows that $T_\xi^{\theta}x\in {\bV^{\geq 1}}$ and thus

\quad (e) {\em
$Q(T_\xi^{[\frac{2n-2m-n_0+2}{2}]}x)=0,\ x\in\bW.$}\\[10pt]
Hence we have

\quad (f) {\em $\lambda_1\leq n+1;$}

\quad (g) {\em
$l_1\leq\max([\frac{2n-2m-n_0+2}{2}],\lambda_1-m),\text{ and } l_1=\lambda_1-m\ \text{ if }[\frac{2n-2m-n_0+2}{2}]\leq\lambda_1-m.$}

Let $\bar{m}$, $\{\bar{v}_i,i\in[0,\bar{m}]\}$, $\{\bar{u}_i,i\in[0,\bar{m}-1]\}$, $\bar{\bW}$, $\bar{\bW}^i$ and $\bar{A}$
be defined for $\bbeta_\xi$ as in subsections \ref{sec-1} and \ref{sec-cen} (if $n=2m$, we choose $\bar{u}_0=u_0^{-2m}$ and then $\bar{u}_{m-1}=u_{m-1}^{-2}$). Note  that we have
$
\bv_{\brm}=v_m;\ \bbeta_\xi(\bar{v}_{\brm-i},\bar{w})=\beta(v_{m-i-1}^{2i+2},\bar{w}), \bar{w}\in\bbW.
$ Hence

\quad (h) {\em $\bar{m}\leq m,\text{ and }\brm=m\text{ iff }v_0^{2m}\neq 0.$}

We first show that\\[5pt]
\indent\quad(i) {\em
$\bW\cap\bV^{\geq-n+1}=\{w\in\bW|T_\xi^{n}w=0,Q(T_\xi^{n-m}w)=0\}$.}\\[5pt]
Let $w\in\bW\cap\bV^{\geq -n+1}$, then $T_\xi^nw\in(\bV^{\geq n+1}+\spn\{v_i,i\in[0,m]\})\cap\bW=\{0\}$. If $n_0>0$, then it follows from (e) that $Q(T_\xi^{n-m}w)=0$.  Assume $n_0=0$. We have $T_\xi^{n-m}w=y+\sum b_iv_i$ for some $y\in\bV^{\geq n+1-2m}$. Since $b_m=\beta_\xi(T_\xi^{n-m}w,u_{m-1})=0$ (note that $\beta_\xi(y,u_{m-1})$=0 as $u_{m-1}\in\bV^{\geq -n+2(m-1)}$), $T_\xi^{n-m}w\in \bV^{\geq 1}$. Hence $Q(T_\xi^{n-m}w)=0$. This shows that the set in the l.h.s of (i) is a subset of that in the r.h.s of (i).

Let $x\in\bW$ be such that $T_\xi^nx=0$ and $Q(T_\xi^{n-m}x)=0$. We show that

\quad (j) {\em$\beta(x^{-n},\bar{A}^n\bbW^{-n})=0.$}\\[10pt]
Suppose $x^{-n}=\sum a_iu_i+\sum b_iv_i+w$, where $w=\pi_\bW(x^{-n})$. Note that if $n_0=0$, then $a_0=\beta(x-x^{-n},v_0)=0$ since $v_0\in\bV^{\geq n}$ and $x-x^{-n}\in\bV^{\geq -n+1}$. It follows that $x+w\in\bV^{\geq -n+1}$ (we use (c)) and thus $T_\xi^nw=0$ and $Q(T_\xi^{n-m}w)=0$. We have $\beta(x^{-n},\bar{A}^n\bbW^{-n})=\bbeta_\xi(x^{-n},\bar{A}^{n-1}\bbW^{-n})=
\beta_\xi(x^{-n},\bar{A}^{n-1}\bbW^{-n})=\beta_\xi(w,\bar{A}^{n-1}\bbW^{-n})$ (in the last equality we use $x^n-w\in\bV^{\geq -n+1}$). It follows from (d) and its proof that
\begin{equation*}\label{eqn-txi}
T_\xi^iw=y_i+\sum_{j\in[0,m]}D_{i+j-m}v_j,
\end{equation*}
for some $y_i\in\bV^{\geq -n+2i}$, where $D_i=0$ for $i\leq[\frac{n}{2}]$, and
$D_i$   for $i>[\frac{n}{2}]$ is the unique number such that $T_\xi^i w+D_iv_m \in \bV^{\geq -n+2i}+\spn\{v_i,i\in[1,m]\}$.
 We show by induction on $i$ that
$$\beta_\xi(w,\bar{A}^{n-1}\bar{\bW}^{-n})=\beta(T_\xi^iw+\sum_{j\in[0,m]} D_{i+j-m}v_j,\bar{A}^{n-i}\bar{\bW}^{-n}).$$
In fact, we have
\begin{eqnarray*} &&\beta_\xi(w,\bar{A}^{n-1}\bar{\bW}^{-n})=\beta(T_\xi w,\bar{A}^{n-1}\bar{\bW}^{-n})=\beta(T_\xi^iw+\sum_{j\in[0,m]} D_{i+j-m}v_j,\bar{A}^{n-i}\bar{\bW}^{-n})\\&&=\bbeta_\xi(T_\xi^iw+\sum_{j\in[0,m]} D_{i+j-m}v_j,\bar{A}^{n-i-1}\bar{\bW}^{-n})=\beta_\xi(T_\xi^{i}w+\sum _{j\in[0,m]}D_{i+j-m}v_j,\bar{A}^{n-i-1}\bar{\bW}^{-n})\\&&
=\beta(T_\xi^{i+1}w+\sum _{{j\in[0,m]}}D_{i+j-m+1}v_j,\bar{A}^{n-i-1}\bar{\bW}^{-n}).
\end{eqnarray*}
It follows that $\beta_\xi(w,\bar{A}^{n-1}\bar{\bW}^{-n})=\beta(T_\xi^nw+
\sum_{j\in[0,m]}D_{n-m+j}v_j,\bar{\bW}^{-n})=0$, since $D_i^2=Q(T_\xi^iw)=0$ for $i\geq n-m+1$ (we use $y_i\in\bV^{\geq 1}$ and (e)); if $n=2m$, then $D_{n-m}=D_m=0$; if $n>2m$, then $D_{n-m}^2=Q(T_\xi^{n-m}w)=0$ (note that $Q(y_{n-m})=0$ as $y_{n-m}\in\bV^{\geq 1}$). Hence (j) holds.

Now assume first that $n=2m$. We have $\bar{m}=m$, $\bar{v}_0=v_0$ (see (h)) and $x^{-2m}\in\bar{\bW}^{-2m}$ (since $\beta(x^{-2m},v_0)=0$). It follows from (j) that $\bar{A}^mx^{-2m}\in \text{Rad}(Q|_{\Ima(\bar{A}^m:\bar{\bW}^{-2m}\to\bar{\bW}^0)})$.  By a similar argument as above one shows that $\beta(\bar{A}^mx^{-2m},\bar{\bW}^0)=\beta((T_\xi^mw)^0,\bar{\bW}^0)$.
Now $(T_\xi^mw)^0\in \bar{\bW}^0$ since $\bbeta_\xi((T_\xi^mw)^0,\bar{u}_{m-1})=\beta_\xi(T_\xi^mw,\bar{u}_{m-1})
=\beta_\xi(T_\xi^mw,u_{m-1})=0$ (note that we have chosen $\bar{u}_0$ such that $\bar{u}_{m-1}=u_{m-1}^{-2}$). Hence $\bar{A}^mx^{-2m}=(T_\xi^mw)^0$ and $Q(\bar{A}^mx^{-2m})=Q((T_\xi^mw)^0)=Q(T_\xi^mw)=0$.
 It follows that $\bar{A}^mx^{-2m}=0$ and thus $x^{-2m}=0$ (see \ref{sec-equiv} (a1)). Hence $x\in\bV^{\geq-n+1}$.

Assume $n>2m$. Then $\bV^{-n}=\bar{\bW}^{-n}$ (note that $n>2m\geq 2\brm$). Since  $\beta(\bar{A}^nx^{-n},\bV^{-n})=\beta(x^{-n},\bar{A}^{n}\bV^{-n})
=0$ (see (j)), $\bar{A}^nx^{-n}=0$ and thus $x^{-n}=0$ (we use that $\bar{A}^n:\bar{\bW}^{-n}\to\bar{\bW}^n$ is an isomorphism). Hence $x\in\bV^{\geq -n+1}$. This complets the proof of (i).

Now we prove the lemma. We have the following cases.

(1) $n=2m$. Then $n_0=0$, $\lambda_1\leq 2m+1$ and $l_1\leq m+1$ (see (f) and (g)); $n=2m\geq m+l_1-1$.
If $l_1\leq m$, then $\lambda_1\leq 2l_1\leq 2m$ and thus for any $x\in \bW$, $T_\xi^nx=0$ and $Q(T_\xi^{n-m}x)=0$; if $l_1=m+1$, then either $m=\lambda_1-l_1=l_1-1$ or $\lambda_1-l_1<m<l_1$. It follows that $H_{\beta_\xi}=\spn\{v_i,i\in[0,m]\}\oplus\spn\{u_i,i\in[1,m-1]\}\oplus\{x\in \bW|T_\xi^nx=0, Q(T_\xi^{n-m}x)=0\}$ and thus $H_{\beta_\xi}\subset\bV^{\geq -n+1}$ (see (c) and (i)). Assume $x=\sum a_iv_i+\sum b_iu_i+\pi_\bW(x)\in\bV^{\geq -n+1}$. We have $b_0=\beta(x,v_0)=0$ since $v_0\in \bV^{\geq n}$.  Then
$\pi_\bW(x)\in\bV^{\geq -n+1}\cap\bW$ (we use (c)) and it follows that $\bV^{\geq -n+1}\subset H_{\beta_\xi}$ (we use (i)).

(2) $n>2m$ and $n_0>0$.  We have $\bV^{-n}=\bar{\bW}^{-n}$ (since $n>2m\geq 2\brm$). Suppose $\lambda_1\leq n$. Then for any $w\in \bW$, $T^n_\xi w=0$ and $Q({T_\xi^{n-m}w})=0$ (see (e)), and it follows that $\bV=\bV^{\geq-n+1}$ (see (c) and (i)) which is a contradiction. Thus $\lambda_1=n+1$ and $l_1=\lambda_1-m$ (see (f) and (g)); $n=m+l_1-1>2m$. Note that  $Q(T_\xi^{l_1-1}\bW)=0$ (see (e)) and thus $\rho=0.$
 Hence  $H_{\beta_\xi}=\spn\{v_i,i\in[0,m]\}\oplus\spn\{u_i,i\in[0,m-1]\}\oplus \ker T_\xi^n$. Now it follows from (c) and (i) that $\bV^{\geq -n+1}=H_{\beta_\xi}$.

 (3) $n>2m$ and $n_0=0$. Let $x=\sum a_iv_i+\sum b_iu_i+\pi_\bW(x)\in\bV^{\geq -n+1}$. We have $b_0=\beta(x,v_0)=0$ (since $v_0\in\bV^{\geq n}$) and thus $\pi_\bW(x)\in\bW\cap\bV^{\geq -n+1}$. It follows that $\bV^{\geq -n+1}=\spn\{v_i,i\in[0,m]\}\oplus\spn\{u_i,i\in[1,m-1]\}\oplus\{x\in \bW|T_\xi^nx=0, Q(T_\xi^{n-m})x=0\}$. We show that $n=m+l_1-1$. Otherwise $l_1\leq n-m$, then $\lambda_1\leq n$ (see (g)) and thus $\bV^{\geq -n+1}=\spn\{v_i,i\in[0,m]\}\oplus\spn\{u_i,i\in[1,m-1]\}\oplus \bW$.  Hence $\dim \bV^{-n}=1$ which is a contradiction ($\dim \bV^{-n}=\dim \bar{\bW}^{-n}$ must be even since $n>2m\geq 2\bar{m}$, see the remark after \ref{sec-equiv} (a2)). Now if $\lambda_1\leq n$, then $\lambda_1-l_1<m<l_1$ and $\bV^{\geq -n+1}=H_{\beta_\xi}$. Assume  $\lambda_1=n+1$, then $m=\lambda_1-l_1<l_1-1$. We show that in this case $\rho\neq 0$. Otherwise, $\bV^n=(\bV^{\geq -n+1})^\p\cap Q^{-1}(0)=\spn\{v_0\}\oplus((\ker T_\xi^{\lambda_1-1})^\p\cap\bW)$ and thus $\dim \bV^n$ is odd which is again a contradiction. It follows that $\bV^{\geq -n+1}=H_{\beta_\xi}$.
This completes the proof of the lemma.

\subsection{}\label{sec-inj}
We prove the injectivity of the map in Proposition \ref{prop-1} by induction on $\dim\bV$.  If $\dim \bV=1$, the statement is clear. Now assume that $\dim\bV\geq 3$. Let $\beta_\xi\in\mathfrak{S}(\bV)$ and let $\bV_*=(\bV^{\geq a})$ and $\tilde{\bV}_*=(\tilde{\bV}^{\geq a})$ be two filtrations in $\mathfrak{F}_o(\bV)$ such that $\beta_\xi\in\eta(\bV_*)$ and $\beta_\xi\in\eta(\tilde{\bV}_*)$. We need to show that $\bV_*=\tilde{\bV}_*$.
If $\beta_\xi=0$, then $\bV^{\geq a}=\tilde{\bV}^{\geq a}=\bV$ for all
$a\leq 0$ and $\bV^{\geq a}=\tilde{\bV}^{\geq a}=0$ for all
$a\geq 1$. Now we assume $\beta_\xi\neq 0$.

Let $\oplus_{a=-n}^n\bV^a$ and $\oplus_{a=-\tilde{n}}^{\tilde{n}}\tilde{\bV}^a$ be $o$-good gradings of $\bV$ such that $\bV^{\geq a}=\oplus_{a'\geq a}\bV^{a'}$ and $\tilde{\bV}^{\geq a}=\oplus_{a'\geq a}\tilde{\bV}^{a'}$. Then  $n=\tilde{n}$ and $\bV^{\geq -n+1}=\tbV^{\geq -n+1}$ (see Lemma \ref{ssec-vn}). Hence $\bV^{n}=\tbV^{ n}$. Let $\bV'=\bV^{\geq -n+1}/\bV^n=\tbV^{\geq -n+1}/\tbV^n$. Then $Q$ induces a nondegenerate quadratic form on $\bV'$. We set $\bV'^{\geq a}=\bV'$ (resp. $\tbV'^{\geq a}=\bV'$) if $a< -n+1$, $\bV'^{\geq a}=\bV^{\geq a}/\bV^n$ (resp. $\tbV'^{\geq a}=\tbV^{\geq a}/\bV^n$) if $a\in[ -n+1,n]$, and $\bV'^{\geq a}=0$ (resp. $\tbV'^{\geq a}=0$) if $a>n$. Then $\bV'_*=(\bV'^{\geq a})\in\mathfrak{F}_o(\bV')$, $\tbV'_*=(\tbV'^{\geq a})\in\mathfrak{F}_o(\bV')$, $\beta_\xi$ induces a symplectic form $\beta_\xi'\in\LS(\bV')$ (see \ref{sec-vu} (i)) and $\beta_\xi'\in\eta(\bV'_*)$, $\beta_\xi'\in\eta(\tbV'_*)$. By induction hypothesis, we have $\bV'^{\geq a}=\bV'^{\geq a}$. It follows that $\bV^{\geq a}=\tbV^{\geq a}$ for $a\geq -n+1$. We have $\bV^{\geq a}=\tbV^{\geq a}=\bV$ for $a< -n+1$. Hence $\bV_*=\tbV_*$.

\subsection{}\label{sec-surj}We prove the surjectivity of the map in Proposition \ref{prop-1} following the arguments used in \cite[2.11]{Lu2}. We can assume that $\tk$ is an algebraic closure of the finite field $\tF_2$ and that $2N+1=\dim\bV\geq 3$. We choose an $\tF_2$ rational structure on $\bV$ such that $Q$ is defined over $\tF_2$. Then the Frobenius map $F$ relative to this $\tF_2$ structure acts naturally and compatibly on $\sqcup_{\bV_*\in\mathfrak{F}_o(\bV)}\eta(\bV_*)$ and $\LS(V)_{nil}$. It is enough to show that for any $n\geq 1$, the map $\Psi_n:(\sqcup_{\bV_*\in\mathfrak{F}_o(\bV)}\eta(\bV_*))^{F^n}\rightarrow\LS(V)_{nil}^{F^n}$, $\beta_\xi\mapsto\beta_\xi$ is a bijection. Since $\Psi_n$ is injective (see \ref{sec-inj}), it suffices to show that $|(\sqcup_{\bV_*\in\mathfrak{F}_o(\bV)}\eta(\bV_*))^{F^n}|=|\LS(V)_{nil}^{F^n}|$. In view of \ref{sec-number} (b),  it is enough to show that

\quad (a) {\em $|(\sqcup_{\bV_*\in\mathfrak{F}_o(\bV)}\eta(\bV_*))^{F^n}|=2^{2nN^2}.$}\\[10pt]
Now the l.h.s of (a) makes sense when $\tk$ is replaced by an algebraic closure of any finite prime field $\tF_{p'}$,  and for $p'\neq 2$, the l.h.s of (a) is equal to ${p'}^{2nN^2}$. Hence it is enough to show that the l.h.s of (a) is a polynomial in $p'^{n}$ with rational coefficients independent of $p'$ and $n$ (following Lusztig, we say that it is {\em universal}).

We now compute the l.h.s of (a)  for general $p'$. A collection of integers $(f_a)_{a\in\mathbb{Z}}$ is called admissible if  $f_{-a}=f_a$, $f_a$ is even for odd $a$, $f_0\geq f_{2}\geq f_{4}\geq\cdots$, $f_{1}\geq f_{3}\geq f_{5}\geq\cdots$ and $\sum_af_a=\dim\bV$. For $(f_a)$ as above, let $\mathcal{Y}_{(f_a)}$ be the set of $\bV_*\in\mathfrak{F}_o(\bV)$ such that $\dim gr_a(\bV_*)=f_a$ for all $a$, where $gr_a(\bV_*)=\bV^{\geq a}/\bV^{\geq a+1}$. We have $|(\sqcup_{\bV_*\in\mathfrak{F}_o(\bV)}\eta(\bV_*))^{F^n}|=\sum_{(f_a)}|\mathcal{Y}_{(f_a)}^{F^n}||\eta(\bV_*)^{F^n}|$, where $\bV_*$ is any fixed element in $\mathcal{Y}_{(f_a)}^{F^n}$. Since $|\mathcal{Y}_{(f_a)}^{F^n}|$ is {\em universal}, it is enough to show that $|\eta(\bV_*)^{F^n}|$ is {\em universal} for any $\bV_*\in\mathcal{Y}_{(f_a)}^{F^n}$. It is easy to see that $|\eta(\bV_*)^{F^n}|=p'^{nd}|(\LS(\bV)_2^0)^{F^n}|$, where $d=\sum_{a<a',a+a'\leq -3}f_af_{a'}+\sum_{a\leq -2}f_a(f_a-1)/2$ is {\em universal} and $\LS(\bV)_2^0$ is defined with respect to an $o$-good grading $\oplus_{a\in\mathbb{Z}}\bV^a$ of $\bV$ such that $\dim\bV^a=f_a$ and $F(\bV^a)=\bV^a$ for all $a$. Let $s'$ be the number of all sequences $U_0\subset U_2\subset U_4\subset\cdots$ of subspaces of $\bV^0$ such that $\dim U_a=f_0-f_a$ and $Q|_{U_a}$ is nondegenerate for all $a$. Let $s''$ be the number of all pairs $(\omega,(U_1,U_3,U_5,\ldots))$ where $\omega$ is a nondegenerate symplectic form on $\bV^{-1}$ and $U_1\subset U_3\subset U_5\subset\cdots$ are subspaces of $\bV^{-1}$ such that $\dim U_a=f_{-1}-f_a$ and $\omega|_{U_a}$ is nondegenerate for all $a$. Let $s_1$ be the number of vector space isomorphisms $\bV^{-1}/U_{2a+1}\to\bV^{2a+1}$ and let $s_2$ be the number of vector space isomorphisms $\bV^{0}/U_{2a}\to\bV^{2a}$. We have that $|(\LS(\bV)_2^0)^{F^n}|=s's''s_1s_2$ is {\em universal}, since $s',s'',s_1,s_2$ are {\em universal} (see \cite[1.2(a),1.2(b)]{Lu3}). This completes the proof of Proposition \ref{prop-1} and thus that of Theorem \ref{thm}.

\subsection{}

Assume that $\tk$ is an algebraic closure of a finite prime field $\tF_p$ and that a split $\tF_p$-rational structure is given on $G$. Then $\Lg$, $\Lg^*$ and $\cN_{\Lg^*}$ have induced $\tF_p$-structures, each $\cO\in\mathfrak{U}_G$ and each  subset $\cN_{\Lg^*}^\cO$ (see \ref{ssec-d1}) are defined over $\tF_p$ (with Frobenius map $F$). As in \cite{Lu2}, it follows from the proof in \ref{sec-surj} that

{\em  for all $n\geq 1$, $|(\cN_{\Lg^*}^\cO)^{F^n}|$ is a polynomial of $p^n$ with integer coefficients independent of $p$ and $n$.}

\section{Examples}
 Let $\bV$ (with $\dim\bV=2N+1$), $Q$, $\beta$ and  $G=SO(\bV)$ be as in subsection \ref{ssec-v}. We fix a good basis  $(e_i)_{i\in[-N,N]}$ of $\bV$ (see \ref{sec-number}). For an $o$-good grading $\bV=\oplus_a\bV^a$ of $\bV$, we denote $\bV_*=(\bV^{\geq a})$, $\bV^{\geq a}=\oplus_{a'\geq a}\bV^{a'}$, the corresponding $Q$-filtration. When $p=2$ and $\beta_\xi\in\LS(\bV)_{nil}$, let $m,\lambda_1,l_1$ and $H_{\beta_\xi}$ be defined for $\beta_\xi$ as in subsection \ref{sec-vu}.

\subsection{}Assume that $N=5$ in this subsection. Let $\bV=\bV^{-4}\oplus\bV^{-2}\oplus\bV^{0}\oplus\bV^{2}\oplus\bV^{4}$ be an $o$-good grading of $\bV$ given by $\bV^{-4}=\text{span}\{e_{-5},e_{-4}\}$, $\bV^{-2}=\text{span}\{e_{-3},e_{-2}\}$, $\bV^{0}=\text{span}\{e_0,e_{-1},e_{1}\}$, $\bV^{2}=\text{span}\{e_{2},e_{3}\}$, $\bV^{4}=\text{span}\{e_{4},e_{5}\}$.

Let $\beta_{\xi_1}\in\LS(\bV)_{nil}$ be such that $\beta_{\xi_1}(e_{i},e_{j})=0$ except $\beta_{\xi_1}(e_{-4},e_{2})=\beta_{\xi_1}(e_{-2},e_{-1})=\beta_{\xi_1}(e_{-2},e_{0})=\beta_{\xi_1}(e_{-5},e_{3})=\beta_{\xi_1}(e_{-3},e_{1})=1$. We have $\beta_{\xi_1}\in\LS(\bV)_2^0\subset\eta(\bV_*)$ (in fact $Ae_0=Ae_{-1}=-e_2$, $Ae_1=-e_3$, $Ae_2=-e_4$, $Ae_3=-e_5$ and $Ae_4=Ae_5=0$). Assume that $p=2$. We have $m=2$, $v_2=e_0,\ v_1=e_2,\ v_0=e_4$; we (can) choose $u_0=e_{-4}$ and then $u_1=e_{-2}$, $\bW=\text{span}\{e_0+e_{-1},e_1,e_{\pm 3},e_{\pm 5}\}$, $e_1\xrightarrow{T_{\xi_1}}e_3\xrightarrow{T_{\xi_1}}e_5\xrightarrow{T_{\xi_1}}0,  e_{-5}\xrightarrow{T_{\xi_1}}e_{-3}\xrightarrow{T_{\xi_1}}e_0+e_{-1}\xrightarrow{T_{\xi_1}}0$; thus $\lambda_1=l_1=3$, moreover $\bV^{\geq -3}=\text{span}\{v_0,v_1,v_2,u_1\}\oplus\{x\in\bW|Q(T_{\xi_1}^2x)=0\}=H_{\beta_{\xi_1}}$.

Let $\beta_{\xi_2}\in\LS(\bV)_{nil}$ be such that $\beta_{\xi_2}(e_{i},e_{j})=0$ except $\beta_{\xi_2}(e_{-4},e_{0})=\beta_{\xi_2}(e_{-4},e_{2})=\beta_{\xi_2}(e_{-5},e_{3})=\beta_{\xi_2}(e_{-3},e_{1})=\beta_{\xi_2}(e_{-2},e_{-1})=1$.  We have $\beta_{\xi_2}\in\eta(\bV_*)$. Assume that $p=2$. Then $m=1$, $v_1=e_0,\ v_0=e_4$, (can choose) $u_0=e_{-4}$, $\bW=\text{span}\{e_0+e_{2},e_{-2},e_{\pm1},e_{\pm3},e_{\pm5}\}$, $e_{-2}\xrightarrow{T_{\xi_2}}e_1\xrightarrow{T_{\xi_2}}e_3\xrightarrow{T_{\xi_2}}e_5\xrightarrow{T_{\xi_2}}0,\ e_{-5}\xrightarrow{T_{\xi_2}}e_{-3}\xrightarrow{T_{\xi_2}}e_{-1}\xrightarrow{T_{\xi_2}}e_0+e_2\xrightarrow{T_{\xi_2}}0$, and thus $\lambda_1=l_1=4$; moreover
 $\bV^{\geq -3}=\text{span}\{v_0,v_1\}\oplus\{x\in\bW|Q(T_{\xi_2}^3x)=0\}=H_{\beta_{\xi_2}}$.

Let $\beta_{\xi_3}\in\LS(\bV)_{nil}$ be such that $\beta_{\xi_3}(e_{i},e_{j})=0$ except $\beta_{\xi_3}(e_{-5},e_{3})=\beta_{\xi_3}(e_{-4},e_2)=\beta_{\xi_3}(e_{-3},e_{1})=\beta_{\xi_3}(e_{-2},e_{-1})=1$. We have $\beta_{\xi_3}\in\LS(\bV)_2^0\subset\eta(\bV_*)$. Assume $p=2$. Then $m=0$, $v_0=e_0$, (can choose) $\bW=\text{span}\{e_{\pm i},i\in[1,5]\}$, $e_{-4}\xrightarrow{T_{\xi_3}}e_{-2}\xrightarrow{T_{\xi_3}}e_1\xrightarrow{T_{\xi_3}}e_3\xrightarrow{T_{\xi_3}}e_5\xrightarrow{T_{\xi_3}}0,\ e_{-5}\xrightarrow{T_{\xi_3}}e_{-3}\xrightarrow{T_{\xi_3}}e_{-1}\xrightarrow{T_{\xi_3}}e_2\xrightarrow{T_{\xi_3}}e_4\xrightarrow{T_{\xi_3}}0$, and thus $\lambda_1=l_1=5$; moreover
 $\bV^{\geq -3}=\text{span}\{v_0\}\oplus\ker T_{\xi_3}^4=H_{\beta_{\xi_3}}$.

 Note that $\beta_{\xi_i},i=1,2,3$ are in the same nilpotent piece; in particular when $p\neq 2$ they are in the same $G$-orbit, but when $p=2$ they are in three distinct  $G$-orbits (see \cite{X2}).

 \subsection{}Assume that $N=8$ and $p=2$ in this subsection.

 Let $\beta_{\xi_1}\in\LS(\bV)_{nil}$ be such that $\beta_{\xi_1}(e_{i},e_{j})=0$ except $\beta_{\xi_1}(e_{-8},e_{0})=\beta_{\xi_1}(e_{-8},e_{1})=\beta_{\xi_1}(e_{-7},e_{4})=\beta_{\xi_1}(e_{-6},e_{3})=\beta_{\xi_1}(e_{-5},e_{2})=\beta_{\xi_1}(e_{-4},e_{-3})=\beta_{\xi_1}(e_{-2},e_{-1})=1$. We have $m=1$, $v_1=e_0,\ v_0=e_8$, (can choose) $u_0=e_{-8}$, $\bW=\text{span}\{e_0+e_{1},e_{-1},e_{\pm i},i\in[2,7]\}$, $e_{-7}\xrightarrow{T_{\xi_1}}e_{-4}\xrightarrow{T_{\xi_1}}e_3\xrightarrow{T_{\xi_1}}e_6\xrightarrow{T_{\xi_1}}0,\ e_{-6}\xrightarrow{T_{\xi_1}}e_{-3}\xrightarrow{T_{\xi_1}}e_{4}\xrightarrow{T_{\xi_1}}e_7\xrightarrow{T_{\xi_1}}0$, $e_{-1}\xrightarrow{T_{\xi_1}}e_{2}\xrightarrow{T_{\xi_1}}e_5\xrightarrow{T_{\xi_1}}0,\ e_{-5}\xrightarrow{T_{\xi_1}}e_{-2}\xrightarrow{T_{\xi_1}}e_0+e_1\xrightarrow{T_{\xi_1}}0$, and thus $\lambda_1=4,\ l_1=3$, $\rho\neq 0$. Let $\bV=\bV^{-3}\oplus\bV^{-1}\oplus\bV^{0}\oplus\bV^{1}\oplus\bV^{3}$ be an $o$-good grading of $\bV$ given by $\bV^{-3}=\text{span}\{e_{i},i\in[-8,-5]\}$, $\bV^{-1}=\text{span}\{e_{i},i\in[-4,-1]\}$, $\bV^{0}=\text{span}\{e_0\}$, $\bV^{1}=\text{span}\{e_{i},i\in[1,4]\}$, $\bV^{3}=\text{span}\{e_{i},i\in[5,8]\}$. It is easy to see that $\beta_{\xi_1}\in\eta(\bV_*)$ and moreover
 $\bV^{\geq -3}=\text{span}\{v_0,v_1\}\oplus\{x\in\bW|T_{\xi_1}^3x=0,\ Q(T_{\xi_1}^2x)=0\}=H_{\beta_{\xi_1}}$.

 Let $\beta_{\xi_2}\in\LS(\bV)_{nil}$ be such that $\beta_{\xi_2}(e_{i},e_{j})=0$ except $\beta_{\xi_2}(e_{-8},e_{3})=\beta_{\xi_2}(e_{-7},e_{2})=\beta_{\xi_2}(e_{-6},e_{0})=\beta_{\xi_2}(e_{-5},e_{1})=\beta_{\xi_2}(e_{-4},e_{-1})=\beta_{\xi_2}(e_{-3},e_{-2})=1$. We have $m=1$, $v_1=e_0,\ v_0=e_6$, (can choose) $u_0=e_{-6}$, $\bW=\text{span}\{e_{\pm i},i\in[1,8]-\{6\}\}$, $e_{-7}\xrightarrow{T_{\xi_2}}e_{-2}\xrightarrow{T_{\xi_2}}e_3\xrightarrow{T_{\xi_2}}e_8\xrightarrow{T_{\xi_2}}0,\ e_{-8}\xrightarrow{T_{\xi_2}}e_{-3}\xrightarrow{T_{\xi_2}}e_{2}\xrightarrow{T_{\xi_2}}e_7\xrightarrow{T_{\xi_2}}0$, $e_{-4}\xrightarrow{T_{\xi_2}}e_{1}\xrightarrow{T_{\xi_2}}e_5\xrightarrow{T_{\xi_2}}0,\ e_{-5}\xrightarrow{T_{\xi_2}}e_{-1}\xrightarrow{T_{\xi_2}}e_4\xrightarrow{T_{\xi_2}}0$, and thus $\lambda_1=4,\ l_1=3$, $\rho=0$.
 Let $\bV=\oplus_{a\in[-3,3]}\bV^{a}$ be an $o$-good grading of $\bV$ given by $\bV^{-3}=\text{span}\{e_{-8},e_{-7}\}$, $\bV^{-2}=\text{span}\{e_{-6},e_{-5},e_{-4}\}$, $\bV^{-1}=\text{span}\{e_{-3},e_{-2}\}$, $\bV^{0}=\text{span}\{e_0,e_{\pm 1}\}$, $\bV^{1}=\text{span}\{e_{2},e_3\}$, $\bV^{2}=\text{span}\{e_4,e_5,e_6\}$, $\bV^{3}=\text{span}\{e_7,e_8\}$.   It is easy to see that $\beta_{\xi_2}\in\LS(\bV)_2^0\subset\eta(\bV_*)$ and  moreover
 $\bV^{\geq -3}=\text{span}\{v_0,v_1,u_0\}\oplus\ker T_{\xi_2}^3=H_{\beta_{\xi_2}}$ (note that  we can choose $w_{**}=0$).

 We have that $\xi_1,\xi_2$ are not in the same nilpotent piece and in particular not in the same $G$-orbit.

\end{document}